 \documentclass[11pt,letterpaper]{amsart}
\usepackage[a4paper,margin=30mm]{geometry}

\usepackage{verbatim,upref,amsxtra,amsthm,xcolor,amssymb,amsxtra}

\usepackage[mathscr]{eucal}
\usepackage{bbm}
\usepackage{dsfont}

\usepackage[initials,nobysame]{amsrefs}

\usepackage{mathtools}
\usepackage[normalem]{ulem}

\usepackage{varioref}

\newcommand{\N}{\mathbb N}

\newcommand{\R}{\mathbb R}

\newcommand{\Dim}{\rm Dim}

\theoremstyle{plain}
\newtheorem{thm}{Theorem}[section]
\newtheorem{lem}[thm]{Lemma}

\newtheorem{cor}[thm]{Corollary}
\newtheorem{exm}[thm]{Example}

\theoremstyle{definition}

\theoremstyle{remark}

\numberwithin{equation}{section}

\date{} 

\begin{document}

\frenchspacing

\setlength\marginparsep{8mm}
\setlength\marginparwidth{20mm}

\title[A topological version of Furstenberg-Kesten theorem]{A topological version of Furstenberg-Kesten theorem} 

\author{Ai Hua Fan}
\address{(A. H. Fan) 
	LAMFA, UMR 7352 CNRS, University of Picardie, 33 rue Saint Leu, 80039 Amiens, France}
\email{ai-hua.fan@u-picardie.fr}


\author{Meng Wu}
\address{(M. Wu)
Department of Mathematical Sciences, P.O. Box 3000, 90014 University of Oulu, Finland
	}
\email{meng.wu@oulu.fi}

	\begin{abstract}
	Let $A(x): =(A_{i, j}(x))$ be a continuous function defined on some subshift of $\Omega:= \{0,1, \cdots, m-1\}^\N$,  taking $d\times d$ non-negative  matrices as values and let $\nu$ be an ergodic $\sigma$-invariant measure
	on the subshift
	where $\sigma$ is the shift map.  Under the condition that $ A(x)A(\sigma x)\cdots A(\sigma^{\ell-1} x)$ is a positive matrix for some 
point  $ x$ in the support of $\nu$ and  some  integer $\ell\ge 1$  and that every entry function $A_{i,j}(\cdot)$ is either identically zero or bounded from below by a positive number which is independent of $i$ and $j$,  it is proved that for any $\nu$-generic point $\omega\in \Omega$,  the limit defining the Lyapunov  exponent 
	$\lim_{n\to \infty} n^{-1} \log \|A(\omega) A(\sigma\omega)\cdots A(\sigma^{n-1}\omega)\|$ exists.  
	\end{abstract}
	
\maketitle

	\section{Introduction}
	Let $\Omega$ be a compact metric space, $T:\Omega \to \Omega$ a continuous map and 
	 $A: \Omega \to {\rm M}_{d} (\mathbb{R})$  a continuous function from $\Omega$ into the space $M_{d} (\mathbb{R})$ of 
	$d\times d$ real matrices. The (maximal)  Lyapunov exponent of $A$ at $\omega\in \Omega$ relative to the dynamics $T$ is defined by the limit 
	\begin{equation}\label{eq:Liap0}
	    L(\omega):=\lim \frac{1}{n} \log \|A(\omega) A(T\omega) \cdots A(T^{n-1}\omega) \|,
	\end{equation}
	if the limit exists.
	For any $T$-invariant  measure $\nu$, Furstenberg and Kesten \cite{FK1960} proved that
	if $A$ takes values in the space ${\rm GL}(d,\mathbb{R})$  of invertible matrices and if $\log^+\|A\| \in L^1(\nu)$, 
	 then the limit defining $L(\omega)$ does exist 
	$\nu$-almost everywhere,  and the limit function $L(\omega)$ is  invariant,  integrable and satisfying
	 \begin{equation}\label{eq:Liap-Max}
	      \Lambda:=\mathbb{E}_\nu  L(\omega) = \lim_n \frac{1}{n}  \mathbb{E}_\nu \log \|A^{(n)}(\omega)\|= \inf_n \frac{1}{n} \mathbb{E}_\nu  \log \|A^{(n)}(\omega)\|
	 \end{equation}   
	 where $ A^{(n)}(\omega) $ denotes the cocycle
	 \begin{equation}\label{eq:cocycle}
	 A^{(n)}(\omega) = A(\omega) A(T\omega) \cdots A(T^{n-1} \omega).
	 \end{equation}
	This result of Furstenberg and Kesten is now a direct consequence of the more general Kingman's subadditive ergodic theorem \cite{Kingman1968}. 
	
	Assume that $\nu$ is ergodic, by the Birkhoff ergodic theorem, $\nu$-almost every point $\omega$ is $\nu$-generic, meaning that
	$$
	      \frac{1}{n}\sum_{k=0}^{n-1} \delta_{T^k \omega}  \rightharpoonup \nu
	$$
	in weak-* topology.  The following question is naturally raised.
	\medskip

	{\bf Question 1.}   Does the limit $L(\omega)$ exist for a given $\nu$-generic point $\omega$ ?
	\medskip
	
	Furstenberg-Kesten theorem does not apply to a fixed individual point $\omega$. 
	In general, the answer to Question 1 is negative, even if the dynamical system $(\Omega, T)$ is uniquely ergodic. There is a counter-example of M. Herman \cite{Herman1981} where $A$ takes values in ${\rm SL}_2(\mathbb{R})$
	and there is another counter-example of P. Walters \cite{Walters1986} where $A$ takes values of $2\times 2$ non-negative matrices (a non-negative matrix $B$ is one whose entries are non-negative real numbers and we write $B\ge 0$.
	By $B>0$ we mean that the entries of $B$ are strictly positive).  Walters' example is constructed on a minimal and uniquely ergodic system whose square is still minimal but not uniquely ergodic. The existence of such system is due to  Veech \cite{Veech1969}. 
	In \cite{Walters1986}, Walters proved the following topological version of Furstenberg-Kesten theorem: 
	if $(\Omega, T)$ is a unique ergodic topological dynamical system with the invariant measure $\nu$ and if $A: X\to {\rm GL}_d(\mathbb{R})$ and $A(\omega) >0$ for all $\omega\in \Omega$, then
	for every $\omega \in \Omega$ the Lyapunov exponent $L(\omega)$ define by (\ref{eq:Liap0}) exists and is equal to $\Lambda$, which is  defined by (\ref{eq:Liap-Max}).   
	One of conditions made by Walters is the unique ergodicity. Another condition is the positivity $A(\omega)>0$, which is crucial and necessary to some extent as the counter-example of Walters shows.
	 A. Furman \cite{Furman1997} found other conditions on  $A(\omega)$ for ensuring the uniformity, i.e. the 
	limit (\ref{eq:Liap0}) exists for every $\omega$ and is uniform in $\omega$.   Furman's conditions are necessary and sufficient when $d=2$. In this direction, Lenz \cite{Lenz2002} obtained a subadditive ergodic theorem on uniquely ergodic subshift systems for the so-called almost additive functions. Lenz  \cite{Lenz2004} also proved that on any subshift $\Omega$ satisfying a condition of  uniform positive weights, then every locally constant
	function $f:\Omega \to {\rm SL} (2, \mathbb{R})$ is uniform. Notice that the  condition of  uniform positive weights implies the minimality and is satisfied by the primitive substitutive systems, even by linearly recurrent systems
	(cf.  \cite{Lenz2004}), which are uniquely ergodic.
	\medskip
	
Without the assumption of  minimality and  unique ergodicity, we  ask the next question:

	{\bf Question 2.}  Under what extra conditions, does the limit $L(\omega)$ exist for a given $\nu$-generic point  $\omega$?
	\medskip
	
	We shall bring a partial answer to Question 2 in the special case where $A$ takes values of non-negative matrices. 
	Let us first state the following 
	result for positive matrix valued function $A(\cdot)$.
	\medskip

	\begin{thm}\label{thm:main1} Let $(\Omega, T)$ be a topological dynamical system and $A: \Omega \to M_{d} (\mathbb{R})$ be a continuous function. 
	Suppose that  \\
	\indent {\rm (i)} \ $A(y) >0$ for all $y\in \Omega$;\\
	\indent {\rm (ii)}  $\omega \in \Omega$ is a $\nu$-generic point for some $T$-invariant measure $\nu$. \\
	Then
	the (maximal)  Lyapunov exponent of $A$ at $\omega$ is well defined by the limit  \eqref{eq:Liap0} and is equal to $\Lambda$ defined by 
	\eqref{eq:Liap-Max}.
	\end{thm}
	
	In \cite{Walters1986}, Walters  stated the conclusion of Theorem \ref{thm:main1}  by making the assumptions that $(\Omega, T)$ is uniquely ergodic and $A$ takes values in ${\rm GL}(d, \mathbb{R})$.
	We observe that both the unique ergodicity  and the invertibility can be dropped  as stated in Theorem \ref{thm:main1}. 
	
	We will prove Theorem \ref{thm:main1} 
	 as corollary of the following topological version of Kingman's theorem.

\begin{thm}\label{thm:main1-1}
Let $(\Omega, T)$ be a topological dynamical system and $\nu$ be a $T$-invariant  measure.  Let $(\varphi_n)_{n\ge 1}$ be a sequence of real valued continuous functions defined on $X$.
Suppose
\begin{itemize}
\item[(i)] $\Lambda:=\lim_{n\to\infty} \frac{\mathbb{E} \varphi_n}{n} \in \mathbb{R}$ exists.
\item[(ii)] $\omega$ is $\nu$-generic.
\item[(iii)] $(\varphi_n)$ is quasi-additive on the orbit of $\omega$, i.e. there exists a sequence of positive numbers $c_n>0$ with $c_n =o(n)$  such that for all integers $n\ge 1, m\ge 1$ and all point $x\in \{T^k\omega: k\ge 0\}$ we have 
\begin{equation}\label{eq:thm1-1-quasi-additive-condition}
 \varphi_n(x)+\varphi_m(T^nx)-c_{n\wedge m}\le \varphi_{n+m}(x)\le \varphi_n(x)+\varphi_m(T^nx)+c_{n\wedge m}. 
\end{equation}
\end{itemize}
Then  the limit $\lim_n \frac{\varphi_n(\omega)}{n}$ exists and is equal to $\Lambda$.
\end{thm}

Here $n\wedge m$ denotes the minimum of $n$ and $m$. Notice that the condition (i) in Theorem \ref{thm:main1-1} is a global property, and it is usually satisfied and is easy to check.
While the condition (iii) is local, meaning that it only concerns the behavior of $\varphi_n$ on the orbit of $\omega$. 
It is possible that the quasi-additivity is satisfied by some generic points but not satisfied by some others. See Example \ref{exm:f(x)} in Section \ref{sect:remarks}.
The proof of Theorem \ref{thm:main1-1} is inspired by that of Walters \cite{Walters1986}.
\medskip

	The positivity condition (i) in Theorem \ref{thm:main1} is too strong. We would like to weaken it. 
	Before stating our main result, 
	let us first recall some basic facts concerning non-negative matrices (cf. \cite{Seneta}). 
	 A nonnegative $d\times d$ matrix 
	 $B=(b_{i, j})\ge 0$ (i.e. the entries $b_{i, j}\ge 0$) is said to be {\em row-allowable} (resp. {\em column-allowable}) if each of its rows (resp. columns) has at least one positive entry.     
 It is said to be {\em allowable} if it is both row-allowable and column-allowable. It is clear that $B\ge 0$ is row-allowable if and only if $Bv>0$ for all positive vector $v>0$. 
 All row-allowable nonnegative matrices form a semigroup, and so do all column-allowable  nonnegative matrices. 
 A sequence $(B_n)$ of non-negative matrices is said to be trivial if $B_1\cdots B_n =0$ for some $n\ge 1$. A sequence of row-allowable (resp. column-allowable) is not trivial.

\medskip

	For the shift dynamics, we will prove a better result than Theorem \ref{thm:main1}. 
	Let $\Omega:=\{0,1, \cdots, m-1\}^\mathbb{N}$.  Recall that  the shift map $\sigma$ on $\Omega$  is defined by $$
	\sigma x = (x_{n+1})_{n\ge 1}\quad {\rm  for}\ \  x=(x_n)_{n\ge 0}.$$ 
	A closed $\sigma$-invariant set $\Sigma$ (i.e. $\sigma(\Sigma)\subset \Sigma$) is called a subshift. The sub-dynamical system 
	$(\Sigma, \sigma)$ is also called a subshift.
	
	The following standard notation will be  used.
 Let $\mathcal{A}=\{0,1, \cdots, m-1\}$, considered as an alphabet.
 We use $$\mathcal{A}^*:=\bigcup_{n=0}^\infty \mathcal{A}^n$$
  to denote  the collection of finite words with letters in $\{0,1, \cdots, m-1\}$
($\mathcal{A}^0$ contains only the empty word).  The length $n$ of a word $u=u_1u_2\cdots u_n \in \mathcal{A}^n$ is denoted by $|u|$.  
Points in $\mathcal{A}^\mathbb{N}$ are considered as
infinite words. Let  $w=(w_k)_{k\ge 1}\in \mathcal{A}^\mathbb{N}$ be an infinite word. For integers $0\le n\le m<\infty$,  we denote
$$
   w_n^m = w_n\cdots w_m.
$$
Such finite words contained in $w$ are also called patterns contained in $w$.
This notation $w_n^m$ can also be used when $w$ is a finite word.  The finites words $w_0^{m-1}$ ($m\ge 1$) are called the prefixes of $w$.
Given a finite word $u\in \mathcal{A}^*$,  $[u]$ denotes the cylinder consisting those $w\in \Omega$ having $u$ as prefix.
For a measure $\mu$ on $\Omega$, ${\rm supp}(\mu)$ stands for its topological support.   
For a matrix $B$, $B_{i,j}$ denote the entries of $B$.
 \begin{thm}\label{Theorem-main-1}
Let $A$ be a continuous function defined on some subshift $\Sigma$ in $\Omega:= \{0,1, \cdots, m-1\}^\N$ taking $d\times d$ non-negative  matrices as values and let $\nu$ be an ergodic measure on $\Sigma$.  Suppose that the following conditions on $A$ and on $\nu$ are satisfied:
\begin{equation}\label{eq:claim1-cond-1}
 \textrm{ there exist }  x\in {\rm supp}(\nu) \textrm{ and } \ell_0\ge 1 \textrm{ such that }  A(x)A(\sigma x)\cdots A(\sigma^{\ell_0-1} x)>0;
\end{equation}
\begin{equation}\label{eq:claim1-cond-2}
 \min_{ 1\le i,j\le d} \min\{A(\omega)_{i,j}: A(\omega)_{i,j}\neq 0,  \omega\in \Omega\} >0. 
\end{equation}
Then for any $\nu$-generic sequence $\omega$, the Lyapunov exponent $L(\omega)$ defined by \eqref{eq:Liap0}  exists and the following alternatives hold:
\begin{itemize}
\item[(i)] there exists $n\in \N$ such that $A(\omega)A(\sigma \omega)\cdots A(\sigma^{n-1}\omega)=0$, so that $L(\omega)=-\infty$;
\item[(ii)] otherwise,   $L(\omega)=\lim_{n\to\infty}\frac{1}{n}\log \mathbb{E}_\nu \|A(x)A(\sigma x)\cdots A(\sigma^{n-1}x)\|\in \mathbb{R}$.
\end{itemize}
\end{thm}

If, furthermore, $A(x)$ is assumed row allowable for all $x$, then the alternative (i) in Theorem \ref{Theorem-main-1} is excluded  and the Lyapunov exponent $L(\omega)$ exits and is finite.
\medskip

We emphasize  the interest  of the existence  of the limit defining $L(\omega)$,  affirmed  in  Theorem \ref{thm:main1} 
and  \ref{Theorem-main-1} for every generic point $\omega$. 
Walters' counter-example shows that  some positivity condition like (\ref{eq:claim1-cond-1}) is necessary to some extent.
The second positivity (\ref{eq:claim1-cond-2}) is also necessary to some extent. Indeed, in Section \ref{sect:remarks}, we shall present examples showing the necessities of the conditions 
(\ref{eq:claim1-cond-1}) and (\ref{eq:claim1-cond-2}), and other remarks.


In the special case where $A$ depends only on the first coordinate, 
the conclusion (ii) in Theorem  \ref{Theorem-main-1} was proved in \cite{Fan2021} (cf. Theorem 1.4 there) but under a stronger condition that   $\omega$ generates a minimal and uniquely ergodic system.
Our proof of Theorem  \ref{Theorem-main-1} will follow the same basic idea as in \cite{Fan2021} by decomposing $\omega$ into return words (see Section \ref{sect:4} for the definition of return word). 
When $\omega$ is minimal, the so-called return words to the cylinder defined by a given word is finite. Without the assumption of minimality, it is possible to have an infinite number of
return words. This creates some difficulty in adapting the arguments in  \cite{Fan2021} to prove Theorem \ref{Theorem-main-1}. However, as we shall see,   we can overcome this
difficult by showing that return words of long lengths are relatively sparse (cf. Lemma \ref{lemma-long-small}).   This is a key point in the proof of Theorem  \ref{Theorem-main-1}.
 
 Let us state a corollary for uniquely ergodic subshifts, which is new as far as we know.

\begin{cor}
Let $Y \subset  \{0,\cdots, m-1\}^{\N}$ is a unique ergodic subshift with  the unique ergodic measure $\nu$ such  that ${\rm supp}(\nu)=Y$.
Let $A_0, A_1, \cdots, A_{m-1}$ be $m (\ge 2)$ non-negative matrices. Suppose that there exists  a word $u=u_0u_1\cdots u_{\ell_0-1}$ ($\ell_0\ge 1$)
such that
$$
     [u]\cap {\rm supp}\, \nu \not=\emptyset, \qquad A_{u_0}A_{u_1}\cdots A_{u_{\ell_0-1}}>0.
$$
Then   the Lyapunov exponent $\lim_{n\to \infty}n^{-1}\log \|A_{\omega_0}A_{\omega_1}\cdots A_{\omega_{n-1}}\|$
exists for every $\omega \in Y$.   
\end{cor}

\medskip

Let us finish this introduction by an application of Theorem \ref{thm:main1} to the multifractal analysis of weighted ergodic averages. Such study was initiated in \cite{Fan1997} and \cite{Fan2021} (See also  \cite{BRS2022} and \cite{BRS2022b}   for some generalizations).
It is one of motivations to the present work. 
Let $T$ be the full shift map on the symbolic space $X=S^\mathbb{N}$ of symbols from $S$ ($q:=|S| \ge 1$), which is equipped with the usual metric, and let $f:X\to \mathbb{R}$ be a function depending
only on the first two coordinates so that we can write $f(x)=f(x_0,x_1)$. For a sequence of weights $w =(w_k)\subset \mathbb{R}$, we
consider the weighted ergodic average
$$
     A^w f(x) = \lim_{n\to \infty} \frac{1}{n}\sum_{k=0}^{n-1} w_k f(T^k x) 
$$
when the limit exists.  Suppose that $(w_k)$ takes values from a finite set $\{v_0, v_1, \cdots, v_{m-1}\}$. Introduce the $q\times q$ positive matrices
$$
            A_j(\beta)  := (e^{\beta v_j f(a, b)})_{(a, b) \in S\times S},  \quad (\forall \beta \in \mathbb{R}, \  \forall j \in \{0, 1, \cdots, m-1\}).
$$
If $(w_k)$ is generic for some shift-invariant measure, by Theorem \ref{thm:main1}, the following limit exists
$$
   \psi(\beta) =\lim_{n\to \infty} \frac{1}{n} \log\|A_{w_1}(\beta)A_{w_2}(\beta) \cdots A_{w_n}(\beta)\| \quad (\forall \beta \in \mathbb{R})
$$ 
and $\psi(\beta)$ is an analytic function of $\beta$ (cf. \cite{Ruelle1979}).  The following corollary follows from Theorem \ref{thm:main1} above and Theorem 1.1 in \cite{Fan2021}. 

\begin{cor}   Suppose that  $(w_k)$ is generic for some shift-invariant ergodic measure. 
If $\alpha =\psi'(\beta)$ for some $\beta \in  \mathbb{R}$, we have
$$
\dim E(\alpha)  = \Dim E(\alpha) = \frac{\psi(\beta) - \alpha\beta}{\log q}
$$
where $E(\alpha) =\{x\in X: A^w f(x)=\alpha \}$, and $\dim E $ and $\Dim E$ denote respectively the Hausdorff dimension and the packing dimension of a set $E$.
 \end{cor}
 
 The assumption that $f$ depends on the first two coordinates is not necessary and we can only assume that $f$ is H\"{o}lder continuous. We can also consider  subshifts of finite type instead of full shift.
 
 We are mainly studying random matrix products determined by a cocycle. But we would like to fix a random point to get a deterministic matrix product. As pointed by Karl Petersen (personal communication), 
 from this point of view,
there is some common point between the present work and the deterministic random walks studied by Aaronson and Keane \cite{AK1982}. 
 \medskip

The rest of the paper is organized as follows. We first prove Theorem \ref{thm:main1-1} in Section \ref{sect:2}, and then Theorem \ref{thm:main1} 
as corollary of 
Theorem 
\ref{thm:main1-1}.  Section  \ref{sect:4}, which is rather long,  is devoted to prove Theorem \ref{Theorem-main-1} and the proofs are different for    continuous or discrete ergodic measure $\nu$. 
In Section \ref{sect:remarks}, we  discuss  the conditions imposed in Theorem \ref{Theorem-main-1} and give other remarks.

\medskip

{\bf Notation.}
For positive numbers $A, B, b, E$ ($b>1$) by the  notation $A=b^{\pm E} B$ we will mean 
$$
               b^{-E}\le \frac{A}{B} \le b^E.
$$
This will simplify estimates on products. Similarly, for real numbers $x, y, c$ ($c>0$), by the  notation $x=y \pm c $ we will mean 
$
                  |x-y| \le c.
$

\section{Proof of Theorem \ref{thm:main1-1}} \label{sect:2}

In this section, we first present a proof of  Theorem \ref{thm:main1-1},  adapting arguments from Walters \cite{Walters1986}.  Then  we shall check that the functions 
$$\varphi_n(x):=\log \|A(x)A(Tx)\cdots A(T^{n-1}x)\|$$   satisfy the condition \eqref{eq:thm1-1-quasi-additive-condition}, which implies Theorem \ref{thm:main1} 
 through Theorem  \ref{thm:main1-1}.

\subsection{Generalized Fekete lemma}
We will need the following generalization of Fekete lemma.

\begin{lem}Let $\{a_n\}_{n\ge 1}$ be a sequence of real numbers such that for all integers $n\ge1, m\ge 1$ we have
\begin{equation}\label{eq:fekete}
   a_{n+m}\le a_n +a_m + c_{n\wedge m}
\end{equation}
where $\{c_n\}_{n\ge 1}$ is a sequence of real numbers such that $\lim_{n\to \infty} n^{-1}c_n =0$. Then the limit  $\lim_{n\to \infty } \frac{a_n}{n} $
exists (it may be $-\infty$).
\end{lem}

\begin{proof} Fix integer $m\ge 1$. For every integer $k\ge 2$ we have 
$$
    a_{km}\le c_m + a_m + a_{(k-1) m}  \le  (k-1) c_m + k a_m.
$$
This estimate remains true when $k=1$.
For any integer $n\ge m$, write $n=km +i$ with $0\le i<m$. We convention that $c_0=0$. Then
$$
    a_n \le a_{km}+ c_i \le \max_{0\le i<k}|c_i| + (k-1) c_m + k a_m.
$$
It follows that
$$
       \frac{a_n}{n}  \le  \frac{1}{k} \frac{\max_{0\le i<m}|c_i|}{m} + \frac{ |c_m|}{m} + \frac{a_m}{m}.
$$
Thus
$$
   \varlimsup_{n\to \infty}   \frac{a_n}{n} \le  \frac{ |c_m|}{m}  + \frac{a_m}{m}.
$$
Now we can conclude by taking liminf on $m$.
\end{proof}

If the inequality in (\ref{eq:fekete}) is reversed, the limit $\lim_{n\to \infty}\frac{a_n}{n}$ exists but it may be $+\infty$. 

Let $a_n = \mathbb{E}_\nu \varphi_n$. If we suppose that the quasi additivity \eqref{eq:thm1-1-quasi-additive-condition} holds for every $x\in \Omega$,  then 
$$
   a_n +a_m -c_{n\wedge m} \le a_{n+m} \le a_n +a_m +c_{n\wedge m}.
$$
Then, by the above Fekete lemma, the following limit exists:
\begin{equation}\label{eq:Lambda}
\Lambda:=\lim_{n\to \infty}\frac{\mathbb{E}_\nu \varphi_n}{n}  \in \mathbb{R}.
\end{equation}
 
\subsection{Proof of Theorem \ref{thm:main1-1}}
We first remark the following consequence of the quasi-additivity condition \eqref{eq:thm1-1-quasi-additive-condition}: Let $p\ge 1$ be an integer.  There exists a constant $c_p^*>0$ such that for 
all $x\in \{T^k\omega: k\ge 0\}$  and all $r\ge 0$ and $0\le i<p$, we have 
\begin{equation}\label{eq:proof-thm1-1-quasi-additive-condition-1}
\varphi_r(T^ix)-c_p^*\le \varphi_{r}(x)\le \varphi_r(T^ix)+c_p^*. 
\end{equation}
It suffices to notice that the condition  \eqref{eq:thm1-1-quasi-additive-condition} implies
\[
\varphi_{i+r}(x)=\varphi_r(x)+\varphi_i(T^rx)\pm c_{i\wedge r}=\varphi_i(x)+\varphi_r(T^ix)\pm c_{i\wedge r},
\]
which implies $\varphi_r(x) - \varphi_r(T^ix) = \varphi_i(x)- \varphi_i(T^r x)\pm 2c_{i\wedge r}$. So,  we can take 
$$ 
c_p^*= 2 \max_{0\le i<p} c_p+ 2 \max_{0\le i <p} \|\varphi_i\|_\infty.
$$

Applying (\ref{eq:proof-thm1-1-quasi-additive-condition-1}) to $x=\omega$,  summing over $0\le i <p$ and then dividing the obtained sum by $p r$, we get 
\begin{equation}\label{eq:mean}
    \left| \frac{1}{r}\, \varphi_r(\omega) -  \frac{1}{p} \sum_{i=0}^{p-1}  \frac{1}{r}\,   \varphi_r(T^i\omega)\right| \le \frac{ c_p^*}{r}.
\end{equation}
Our objective is to prove that  $r^{-1}\varphi_r(\omega)$ tends to $\Lambda$ which is defined by (\ref{eq:Lambda}). By (\ref{eq:mean}),  we are led to study the arithmetic mean of the $p$ quantities $\frac{1}{r}\,   \varphi_r(T^i\omega)$ for $0\le i <p$, which will approach to $\Lambda$, as we shall prove. 
Let us first look at the arithmetic mean in (\ref{eq:mean}) for $r=np$. 

For $\epsilon >0$, take a large $p$ such that 
$c_p<p\epsilon$ and 
\begin{equation}\label{eq:Elogg}
   \left| \frac{1}{p} \mathbb{E}_{\nu}  \varphi_p(x) - \Lambda \right| <\epsilon.
\end{equation}
Let  $g(\cdot) =\varphi_p(\cdot)$.
Now,  by \eqref{eq:thm1-1-quasi-additive-condition} which is used $n-1$ times, for every point $x$ on the orbit of $\omega$ we get
$$
               -(n-1)c_p+ \sum_{j=0}^{n-1} g(T^{pj} x) \le   \varphi_{np}(x) \le \sum_{j=0}^{n-1} g(T^{pj}x)+(n-1)c_p.
$$
Replacing $x$ by $T^i \omega$,  summing the above inequalities over $0\le i<p$ and then dividing by $n p^2$, we get 
\begin{equation}\label{eq:1-epsilon}
                                    \left|  \frac{1}{p}\sum_{i=0}^{p-1} \frac{1}{np}   \varphi_{np}(T^i \omega) -   \frac{1}{p}\cdot \frac{1}{np}\sum_{t=0}^{np-1} g(T^t \omega)  \right| \le   \frac{(n-1)c_p}{pn}  
                                    <\epsilon.
\end{equation}
where we have used the fact that $p$ was chosen so that $c_p<p\epsilon$. Since $\omega$ is generic, we have 
$$
\lim_{n\to \infty }\frac{1}{np}\sum_{t=0}^{np-1} g(T^t \omega)  = \mathbb{E}_\nu g.
$$
This, together with  (\ref{eq:Elogg}), implies that there exists  an integer $N_p$ such  that for $n\ge N_p$, we have
\begin{equation}\label{eq:2-epsilon}
      \left|  \frac{1}{p}\cdot \frac{1}{np}\sum_{t=0}^{np-1}  g(T^t \omega) - \Lambda \right|< 2 \epsilon.
\end{equation}
From (\ref{eq:1-epsilon}) and  (\ref{eq:2-epsilon}) we get immediately 
 \begin{equation}\label{eq:mean2}
                    \forall n\ge N_p, \quad                \left|  \frac{1}{p}\sum_{i=0}^{p-1} \frac{1}{np}\varphi_{np}(T^i \omega) -   \Lambda  \right| \le   3\epsilon.
\end{equation}    

For a general $r$,  we write $r=np + q$ with $0\le q <p$.  We have,  again by \eqref{eq:thm1-1-quasi-additive-condition},
$$  
     \varphi_r(x) =\varphi_{np}(x) +\varphi_q(T^{np}x)\pm c_q=\varphi_{np}(x) +O(1),
$$
where O(1) is a quantity  uniformly bounded  in $x$ and $r$, but depending on $p$. Hence, applying this to $x=T^i \omega$ for $0\le i<p$, we can deduce
$$
    \frac{1}{np^2}\sum_{i=0}^{p-1}\varphi_r(T^i \omega)  = \frac{1}{np^2}  \sum_{i=0}^{p-1}\varphi_{np}(T^i \omega) +  O((np)^{-1}).
$$
 This, together with (\ref{eq:mean2}), allows us to take a larger $N_p$ if necessary such that  
 $$
    \forall n \ge N_p, \qquad  \left| \frac{1}{np^2}\sum_{i=0}^{p-1}\varphi_r(T^i \omega)  - \Lambda \right| <4\epsilon.
$$
Since $\frac{np}{r} \to 1$, we can even claim that 
 $$
    \forall n \ge N_p, \qquad  \left| \frac{1}{p}\sum_{i=0}^{p-1} \frac{1}{r}\varphi_r(T^i \omega)  - \Lambda \right| <5\epsilon.
$$
This, together with (\ref{eq:mean}), finishes the proof of $\lim \frac{1}{r} \varphi_r(\omega) = \Lambda$.


\subsection{Two elementary lemmas on non-negative matrices}
Before giving the proofs of Theorem \ref{thm:main1},  
let us recall the following elementary lemmas which are fundamental for our subsequent studies on the products of non-negative matrices.  It  will also be  used to prove 
Theorem \ref{Theorem-main-1}.

Recall that for non-negative matrix $B=(b_{i,j})$, we use the matrix norm $\|B\|$ defined by
  $$
 \|B\|=\sum_{i, j} b_{i, j} = \,^t\mathbf{1} B \mathbf{1},$$
  where $\mathbf{1} $ is the column vector having $1$ as entries and $^t\mathbf{1}$ denotes the transpose of $\mathbf{1}$.
 
 \begin{lem} \label{lem:elem}Let $P=(p_{i,j})>0$ be a positive $d\times d$ matrix. There exists a positive constant $0<c(P)\le 1$ depending on $P$ such that for all non-negative $d\times d$ matrices  $L$ and $R$ we have
 $$
    c(P) \|L\| \|PR\| \le  \|LPR\|\le \|L\| \|PR\|.
 $$
 We can take $c(P) = d^{-1}\frac{\min_{i, j} p_{i, j}}{\max_{i, j} p_{i, j}}$.
 \end{lem}
 \begin{proof} 
 Let $\widetilde{P} = M P$ with $M=\mathbf{1} \,^t\mathbf{1}$ (the matrix having $1$ as entries). Notice that 
 $$
       \widetilde{p}_{i, j} = \sum_{k} p_{k,j} \le d \max_{m, n} p_{m, n} \le  d \frac{\max_{m, n} p_{m, n}}{\min_{m, n} p_{m, n}}
 p_{i, j} = \frac{p_{i, j}}{c(P)}. $$
 In other words, $P\ge c(P) \widetilde{P}$. Then
 $$
   \|LPR\| \ge c(P) \|L\widetilde{P} R\| =c(P)\cdot\, ^t\!\mathbf{1} L \mathbf{1} \cdot \,^t\!\mathbf{1} P R\mathbf{1} = c(P) \|L\| \|PR\|.
 $$
 The other inequality can be checked  by the definition of the norm.
 \end{proof}

 \begin{lem}\label{lem:base} Let $B_1, B_2, \cdots, B_n$ be $d\times d$ non-negative matrices such that $B_1B_2\cdots B_n\not=0$. Then 
$$
    -c n \le \log \|B_1B_2\cdots B_n\|\le c n
$$
where $c= \max \{|\log a_*|, |\log (a^* d^2)|\}$ with $a_*$ (resp. $a^*$) being the minimum (resp. maximum) of all non-zero entries of all $B_1, B_2, \cdots, B_n$.
\end{lem}

\begin{proof} 
There is   at least one entry of $B_1B_2\cdots B_n$ which is non-zero. So $$
\|B_1B_2\cdots B_n\|\ge a_*^{n}.$$ 
On the other hand,
we  have 
\[\|B_1B_2\cdots B_n\| \le \|(a^* \mathbf{1} \,^t \mathbf{1})^{n}\| \le (a^* d^2)^{n}.\]
We can then conclude. 
\end{proof}

\subsection{Proof of Theorem \ref{thm:main1}}


 Let $\varphi_n(x)=\log \|A^{(n)}(x)\|$. 
 By Theorem  \ref{thm:main1-1}, we only need to check  that the sequence $\{\varphi_n(x)\}$ satisfies the quasi-additivity condition \eqref{eq:thm1-1-quasi-additive-condition}
 for every point $x\in \Omega$.   Indeed, since $A(T^nx)$ is positive,  by Lemma  \ref{lem:elem},  we have
 \[
c\|A^{(n)}(x)\| \|A^{(m)}(T^nx)\|\le   \|A^{(n+m)}(x)\|\le \|A^{(n)}(x)\| \|A^{(m)}(T^nx)\|
 \]
for some constant $0<c\le 1$ depending on the function $A(\cdot)$. 
 Taking logarithm,     we get the quasi-additivity \eqref{eq:thm1-1-quasi-additive-condition}  for $\{\varphi_n(x)\}$ defined above,  with $c_n=|\log c|$.

\section{ Proofs of Theorem \ref{Theorem-main-1}  }\label{sect:4}

It will be convenient to  introduce the following notation.  For integers $0\le n< m$ and an arbitrary point $x\in \Omega$, denote
\[A^{(n,m)}(x):= A(\sigma^{n}x)\cdots A(\sigma^{m-1}x) = A^{(m-n)}(T^nx).  \]
This notation $A^{(n,m)}(x)$ generalizes $A^{(m)}(x)$, because
 $A^{(m)}(x)=A^{(0,m)}(x)$. 

When there exists $n\in \N$ such that $A^{(n)}(\omega)=0$,  we trivially have $L(\omega)=-\infty$.  Thus in the following we shall  suppose that $A^{(n)}(\omega)\neq 0$ for all $n\in \N$.  Then we shall show that 
under the conditions \eqref{eq:claim1-cond-1} and \eqref{eq:claim1-cond-2}, the Lyapunov exponent $L(\omega)$ defined by \eqref{eq:Liap0}  exists and the  alternative (ii) of the conclusion of  Theorem \ref{Theorem-main-1} holds.


\subsection{Positivity on a cylinder} 
\label{sect4-1}
Our proof starts with a trivial  consequence  of the condition \eqref{eq:claim1-cond-1}, that we state as a lemma.  
\begin{lem}\label{lemma-implication-condition-1}
There exists a cylinder $[u]$ determined by a word $u$ with  $|u|\ge \ell_0$, for which we have 
\begin{equation}\label{eq:implication-condition-1}
 \nu([u])>0 \textrm{ and } b:=\min_{1\le i,j\le d} \min_{x\in [u]}A^{(\ell_0)}(x)_{i,j}>0.
\end{equation}
\end{lem}
\begin{proof}
Since the function $z\mapsto A^{(\ell_0)}(z)  $ is continuous, so are the entry functions $A^{(\ell_0)}(z)_{i, j}$.  The  condition \eqref{eq:claim1-cond-1} means that 
$A^{(\ell_0)}(x)_{i, j} >0$ for all $i$ and $j$ at a point $x\in {\rm supp}\, \nu$.  So, there is an neighborhood of $x$, say $[u]$, on which all the functions $A^{(\ell_0)}(\cdot)_{i, j}$ are positive.   
The least of the minimal values on $[u]$ of $A^{(\ell_0)}(\cdot)_{i, j}$ is positive. If we like, we can assume $|u|\ge \ell_0$. Since $x\in {\rm supp}\, \nu$, we have $\nu([u])>0$.
\end{proof}

In the following subsections $\S$\ref{subsect:decomp}-$\S$\ref{subsection-proof-lem-long-small}, we first give a proof for  Theorem \ref{Theorem-main-1}  when  the measure $\nu$ is non-atomic. The proof for  atomic measure $\nu$ is easier and will be given later in the subsection $\S$\ref{sect:atomic}.   Our proof of Theorem \ref{Theorem-main-1}  consists of a series of lemmas and is based on the decomposition of $\omega$ into return words with respect to  some  sub-cylinder sets
of the  cylinder set $[u]$ which appears in Lemma \ref{eq:implication-condition-1}.

\subsection{Decomposition of $\omega$ into return words}\label{subsect:decomp}

We start with the following lemma, which asserts that $\sigma^{|u|}$ pushes ${\rm supp}(\nu)\cap [u]$ out of $[u]$. 

\begin{lem} \label{lem:1}There  exists a point $z\in {\rm supp}(\nu)\cap [u]$ such that $\sigma^{n_0} z\notin [u]$, where $n_0= |u|$.
\end{lem}

\begin{proof} We prove the assertion by contradiction. 
Suppose  $\sigma^{n_0}({\rm supp}(\nu) \cap [u])\subset[u]$.  Since  the support of the invariant measure $\nu$ is $\sigma$-invariant, i.e. $\sigma( {\rm supp}(\nu)) \subset  {\rm supp}(\nu)$,  
we  get  $$\sigma^{n_0}({\rm supp}(\nu)\cap [u])\subset  {\rm supp}(\nu) \cap [u],$$  which implies by iteration  that $ {\rm supp}(\nu) \cap [u]=\{u^\infty\}$.  This contradicts the assumption that $\nu$ is non-atomic.
\end{proof}

Let us now fix $z\in [u]\cap {\rm supp}(\nu)$ such that $\sigma^{n_0} z\notin [u]$.  Such a point $z$ exists by Lemma \ref{lem:1}. Since $\nu$ is non-atomic, we have 
$$\lim_{k\to\infty}\nu([z_0^k])=0.$$
In the following we fix a large $k_0$  and take $v=z_0^{k_0}$. Notice that the measure $\nu([v])$ can be as small as we want if $k_0$ is large enough. 
Actually we will need a sequence of $k_0$ tending to the infinity.

Now let us consider the set of the {\em return times} of $\omega$ into the cylinder $[v]$:
 \[
 T(v):=\{k\ge 1: \sigma^{k}\omega\in [v]\}=\{\tau_i\}_{i\ge 0}, \quad {\rm with} \ \  1\le \tau_0<\tau_1<\cdots.
 \]
 Since $\omega$ is $\nu$-generic and $\nu([v])>0$, the set $T(v)$ is infinite.  So, we can decompose the sequence $\omega$ as
\begin{equation}\label{eq:decomp}
\omega=\zeta_0 \zeta_1\zeta_2\cdots\zeta_i \cdots, \quad {\rm with} \ \ \zeta_0=\omega_0^{\tau_0-1}, \zeta_i = \omega_{\tau_{i-1}}^{\tau_i-1}\  (\forall i \ge 1). 
\end{equation}
The collection of {\em return words} of $\omega$ with respect to $v$ is defined to be
\begin{equation}\label{eq:return-words-1}
R(v)=\left\{ \zeta_j: j\ge 1\right\}.
\end{equation}
Notice that for any $j\ge 1$, $\sigma^{\tau_{j-1}} \omega$ (i.e. $\zeta_j\zeta_{j+1}\cdots$) has $v$ as prefix. 
It is possible that $R(v)$ is finite and $\tau_{k+1}-\tau_k$ is bounded. That is the case when $\omega$ is minimal.  
It is also possible that $R(v)$ is  infinite and then there are return words  of lengths tending to infinity. This creates some difficulties for us to prove Theorem \ref{Theorem-main-1}. 
However these difficulties can be overcome, because it can be proved that there are long return words but not too much (cf. Lemma  \ref{lemma-long-small}). 
\medskip

Firstly, Lemma \ref{lem-exist-frequency-1} below  shows that roughly speaking, when we shift $\omega$, we see $v$ in $\omega$ every $\frac{1}{\nu([v])}$ times, because
$$
    \lim_{i\to\infty} \frac{\tau_i}{i} = \lim_{i\to\infty} \frac{1}{i}\sum_{k=0}^{i-1} (\tau_{k+1}-\tau_k)
$$ 
and so the interval between two consecutive return times is $\frac{1}{\nu([v])}$ in mean.
Or, in other words, most of return words are of bounded length. Secondly,
Lemma \ref{lemma-long-small} below shows that there are  few  long return words to some extent.


\begin{lem}\label{lem-exist-frequency-1}
We have 
\begin{equation}\label{eq:lem-exist-frequency-1}
\lim_{i\to\infty}\frac{i}{\tau_i}=\nu([v]).
\end{equation}
Consequently, 
\begin{equation}\label{eq:lem-exist-frequency-2}
\lim_{i\to\infty}\frac{\tau_{i+1}}{\tau_i}=1.
\end{equation}
\end{lem}

\begin{proof}
Notice that the number $i$ appearing in the decomposition (\ref{eq:decomp}) is the number of return times to the cylinder $[v]$ along the orbit $\{\sigma^j\omega\}_{1\le j\le \tau_i}$. 
Thus
\[
\frac{i}{\tau_i}=\frac{1}{\tau_i}\sum_{j=1}^{\tau_i}{\bf 1}_{[v]}(\sigma^j\omega),
\]
which tends to $\nu([v])$, by the $\mu$-genericity of $\omega$. The limit 
 \eqref{eq:lem-exist-frequency-1} is thus proved. (\ref{eq:lem-exist-frequency-2}) follows from  (\ref{eq:lem-exist-frequency-1})  because $\nu([v]) >0$.
\end{proof}

For any integer $M\ge 1$, consider the finite set of return words having their lengths not exceeding $M$:
$$F^M_v=\{\zeta\in R(v): |\zeta|\le M\}.$$
Lemma \ref{lemma-long-small} below shows that among the first $i$ return words, those having long lengths have their total length negligible    with respect to $\tau_i$.
When $R(v)$ is a finite set (that is the case when $\omega$ is minimal),  Lemma \ref{lemma-long-small} is trivial  because $R(v) \setminus F_v^M$ is empty for large $M$.

\begin{lem}\label{lemma-long-small}
	We have 
	$$\lim_{M\to\infty}\varlimsup_{i\to\infty}\frac{1}{\tau_i}\sum_{j=1}^{i} \left |\zeta_j\right|{\bf 1}_{R(v)\setminus F_v^M}(\zeta_j)=0.$$
\end{lem}
The proof of Lemma  \ref{lemma-long-small} is technical. We postpone its proof in Section \ref{subsection-proof-lem-long-small}. 
\medskip

Let us make the following remark about the infinity of return words for a generic point $\omega$ relative to a Bernoulli measure $\nu$, say $\nu([j])=\frac{1}{m}$ for $0\le j <m$. 
Notice that  $\omega$ contains all possible finite words. Assume a word $v$ starts with the letter $1$.  For any integer $t\ge 1$, as $0^t$ (the word of length $t$ consisting of $0$'s) appears infinitely many times in 
$\omega$,  there are infinitely many return times $\tau_i \ge t$ (for returning to $[v]$). If the word $v$ starts with the letter $j$, it suffices to consider the occurrences of $k^t$ with $k\not=j$.

\subsection{Quasi-multiplicativity of $\|A^{(n)}(\omega) \|$}\label{sect4-3}


With the help of Lemma \ref{lem:elem}, the condition (\ref{eq:claim1-cond-1}) implies the following quasi-multiplicativity of $\|A^{(n)}(\omega)$ along return words.
But in general, it is not possible to compare $\|A^{(n+m)}(\omega)\|$ with $\|A^{(n)}(\omega)\| \|A^{(n, n+m)}(\omega)\|$  for all $n$ and $m$.

\begin{lem}\label{lem-quasi-multiplicativity-1}
There exists $0<c_1\le 1$,  depending only on $A$ and $u$,   such that for every $j\ge 1$ and every $\ell\ge |u|$, we have
 \begin{equation}\label{eq:lem-quasi-multiplicativity-0}
c_1\le \frac{\|A^{(\tau_j+\ell)}(\omega)\|}{  \|A^{(\tau_j)}(\omega)\|  \|A^{(\tau_j,\tau_j+\ell)}(\omega)  \|}\le 1.
\end{equation}
Consequently, for each $i\ge 1$, we have 
 \begin{equation}\label{eq:lem-quasi-multiplicativity-1}
c_1^i \cdot \prod_{j=0}^{i}\|A^{(\tau_{j-1},\tau_j)}(\omega) \|\le  \|A^{(\tau_i)}(\omega)\|\le \prod_{j=0}^{i} \|A^{(\tau_{j-1},\tau_j)}(\omega)\|,
\end{equation}
where we make the convention $\tau_{-1}=0$.
\end{lem}

\begin{proof}
By the definition of $\tau_j$,  $\sigma^{\tau_j}(\omega)\in [v]\subset  [u]$.  The assumption  $\ell\ge |u|$ then implies that $u$  is a prefix of the word 
$\sigma^{\tau_j}(\omega)|_0^{\ell-1}$ so that $\sigma^{\tau_j}(\omega) \in [u]$.  Recall  that $|u|\ge \ell_0$. Hence, by Lemma  \ref{lemma-implication-condition-1},  all the entries of the matrix $A^{(\tau_j, \tau_{j}+\ell_0)}(\omega)$ are bounded from below by some positive number $b>0$.  Now observe that the matrix  $A^{(\tau_j+\ell)}(\omega)$ can be read as  $LPR$ with 
  $$
L=A^{(\tau_j)}(\omega), \quad P=A^{(\tau_j, \tau_{j}+\ell_0)}(\omega), R= A^{(\tau_j+\ell_0, \tau_{j}+\ell)}(\omega).
$$ Applying Lemma \ref{lem:elem} we get
the  inequalities in (\ref{eq:lem-quasi-multiplicativity-0}).

Prove now \eqref{eq:lem-quasi-multiplicativity-1}.   We only need to prove the first inequality. Recall that $v=z_0^{k_0}$ (with large $k_0$) and $\sigma^{n_0}z\notin [u]$.    This implies that  we must have $\tau_{j+1}- \tau_j\ge n_0$ for all $j\ge 1$. 
By \eqref{eq:lem-quasi-multiplicativity-0}, we first get 
$$c_1 \|A^{(\tau_{i-1})}(\omega) \| \|A^{(\tau_{i-1},\tau_i)}(\omega) \|\le  \|A^{(\tau_i)}(\omega)\|.
$$
Then, inductively, we get the first inequality in (\ref{eq:lem-quasi-multiplicativity-1}).
\end{proof}

Note the useful fact that the constant $0<c_1\le 1$ above only depends on $A$ and $u$, but  not on $v$.  In fact, we can take $$
c_1=\min_{x\in [u]}c(A^{(\ell_0)}(x))$$
where  $c(\cdot)$ is the function appearing in Lemma \ref{lem:elem}.

The condition (\ref{eq:claim1-cond-2}) implies the following rough estimation of $\log \|A^{(n, m)}\|(\omega)$, which is nothing but a restatement of lemma \ref{lem:base}.

\begin{lem}\label{lemma-basic-norm-estimate-1}
There exists $0<c_2<\infty$ (depending only on $A$ and $\omega$) such that for all $0\le n\le m$, we have
$$-c_2(m-n) \le \log \|A^{(n,m)}(\omega)\| \le c_2 (m-n).$$
\end{lem}

\subsection{Proof of Theorem \ref{Theorem-main-1} when $\nu$ is not atomic} \label{sect4-4}
\ \,

{\em Step 1. It suffices to prove the existence of $\lim_{i\to \infty}\frac{\log \|A^{(\tau_i)}(\omega)\|}{\tau_i}$.}
Assume $n\ge \tau_1+|u|$. There exists a unique $i\ge 1$ such that 
\[\tau_i+|u|\le n<\tau_{i+1}+|u|.\]
By the estimates  \eqref{eq:lem-quasi-multiplicativity-0} in Lemma \ref{lem-quasi-multiplicativity-1} and Lemma  \ref{lemma-basic-norm-estimate-1},  we have
\[
\log \|A^{(\tau_i)}(\omega)\|-(n-\tau_i)c_2 \le \log \|A^{(n)}(\omega)\| \le  \log \|A^{(\tau_i)}(\omega)\|+(n-\tau_i)c_2.
\]
By Lemma \ref{lem-exist-frequency-1}, we have $\frac{\tau_{i+1}}{\tau_i}\to 1$, which implies $\frac{n}{\tau_i} \to 1$. So, we have 
\begin{equation}\label{eq:1}
 \varliminf_{i \to\infty} \frac{\log \|A^{(\tau_i)}(\omega)\|}{\tau_i} \le \varliminf_{n\to\infty} \frac{\log \|A^{(n)}(\omega)\| }{n}  \le 
\varlimsup_{n\to\infty} \frac{\log \|A^{(n)}(\omega)\| }{n}  
\le \varlimsup_{i \to\infty} \frac{\log \|A^{(\tau_i)}(\omega)\|}{\tau_i}.
\end{equation}

{\em Step 2. Conversion to the existence of $ \lim_{i\to \infty}\frac{1}{\tau_i} \sum_{j=0}^{i}\log \|A^{(\tau_{j-1},\tau_j)}(\omega)\|$.}
The quasi-multiplicativity \eqref{eq:lem-quasi-multiplicativity-1}  implies
\begin{equation}\label{eq:proof-main-thm-qm-1}
\frac{\log \|A^{(\tau_i)}(\omega)\|}{\tau_i} =
\frac{1}{\tau_i} \sum_{j=0}^{i}\log \|A^{(\tau_{j-1},\tau_j)}(\omega)\|+O\left(\frac{i}{\tau_i}\right),
\end{equation}
where the constant involved  in $O\left(\frac{i}{\tau_i}\right)$ is $|\log c_1|$ and is  independent of $v$.  
 Then we can make $O\left(\frac{i}{\tau_i}\right)$ as small as we want by choosing $v$ with long length. Indeed, as $v=z_0^{k_0}$ and  $\lim_{k\to\infty}\nu([z_1^k])=0$,    for any arbitrarily  small $\epsilon>0$, we  can take a sufficiently large $k_0\gg 1$ such that $\nu[z_0^{k_0}]<\frac{\epsilon}{2}$. 
 Then \eqref{eq:lem-exist-frequency-1} implies that
\begin{equation}\label{eq:2}
\frac{i}{\tau_i}\le \epsilon \textrm{ if } i\gg 1.
\end{equation}
If we can prove the existence of the following limit
\begin{equation}\label{eq:limit}
    \lim_{i\to \infty}\frac{1}{\tau_i} \sum_{j=0}^{i}\log \|A^{(\tau_{j-1},\tau_j)}(\omega)\|,
\end{equation}
then from from (\ref{eq:1}), (\ref{eq:proof-main-thm-qm-1}), (\ref{eq:2}) and (\ref{eq:limit}), we get 
$$
    \varlimsup_{n\to\infty} \frac{\log \|A^{(n)}(\omega)\| }{n} - \varliminf_{n\to\infty} \frac{\log \|A^{(n)}(\omega)\| }{n}  =O(\epsilon).
$$
Since $\epsilon >0$ is arbitrary, we finish the proof. 

{\em Step 3. The proof of the existence of the limit in \eqref{eq:limit}.}
Recall that $R(v)=\{\zeta_1, \zeta_2, \cdots\}$ denotes  the set  of  return words to $[v]$, defined by \eqref{eq:return-words-1}.
See the decomposition (\ref{eq:decomp}) of $\omega$ into $\zeta_j$'s.
For any fixed return word $\zeta \in R(v)$,  let us pay attention to those $0\le j\le i$  such that  $\zeta_j=\zeta$ in order to get
\[
\sum_{j=0}^{i}{\bf 1}_{\zeta}(\zeta_j)\log \|A^{(\tau_{j-1},\tau_j)}(\omega)\|= \sum_{\ell=0}^{\tau_i}{\bf 1}_{[\zeta]}(\sigma^\ell\omega )\log \|A^{(|\zeta|)}(\sigma^\ell\omega)\|.
\]
Then, by the $\mu$-genericity of $\omega$, we have
\[
\lim_{i\to\infty}\frac{1}{\tau_i}\sum_{j=0}^{i}{\bf 1}_{\zeta}\left(\zeta_j\right)\log \|A^{(\tau_{j-1},\tau_j)}(\omega)\|
=\int_{[\zeta]}\log \|A^{(|\zeta|)}(x)\|  d\nu(x).
\]
The last integral is a finite real number because the integrand is bounded according to Lemma \ref{lemma-basic-norm-estimate-1}.
It follows that for any integer $M\ge 1$ we have
\begin{equation}\label{eq:limit-approx-formula-1}
\lim_{i\to\infty}\frac{1}{\tau_i}\sum_{j=0}^{i} {\bf 1}_{F_v^M}\left(\zeta_j\right) \log \|A^{(\tau_{j-1},\tau_j)}(\omega) \|=\sum_{\zeta\in F_v^M}\int_{[\zeta]}\log \|A^{(|\zeta|)}(x)\|  d\nu(x),
\end{equation}
where $F_v^M=\{\zeta \in R(v): |\zeta|\le M\}$, which is a finite set.
On the other hand, by  Lemma \ref{lemma-basic-norm-estimate-1} and Lemma \ref{lemma-long-small},  for any
arbitrary small  $\epsilon>0$, there exists $M\ge 1$ such that 
 \begin{equation}\label{eq:limsup-bound-1}
 \limsup_{i\to\infty}\frac{1}{\tau_i}\sum_{j=0}^{i}  {\bf 1}_{R(v)\setminus F_v^M}\left(\zeta_j\right) \left|\log  \|A^{(\tau_{j-1},\tau_j)}(\omega) \| \right|\le \epsilon.
 \end{equation}
 From (\ref{eq:limit-approx-formula-1}) and (\ref{eq:limsup-bound-1}), we get 
 \begin{equation}\label{eq:3}
    \varlimsup_{i\to\infty}\frac{1}{\tau_i}\sum_{j=0}^{i}  \log \|A^{(\tau_{j-1},\tau_j)}(\omega) \|-  \varliminf_{i\to\infty}\frac{1}{\tau_i}\sum_{j=0}^{i}  \log \|A^{(\tau_{j-1},\tau_j)}(\omega) \| \le \epsilon
 \end{equation}
 Since $\epsilon>0$ is arbitrary, we have thus proved the existence of the limit in (\ref{eq:limit}).

{\em Step 4. The identification of the limit.}
We have actually proved that  for any $\nu$-generic point $\omega$ with $A^{(n)}(\omega)\neq 0$ for all $n\in \N$, the Lyapunov exponent $L(\omega)$ exists and equals to 
\begin{equation}\label{eq:Liap-Formula}
   L(\omega) =\lim_{k_0\to \infty} \lim_{M\to \infty}\sum_{\zeta\in F_v^M}\int_{[\zeta]}\log \|A^{(|\zeta|)}(x)\|  d\nu(x),
\end{equation}
where $k_0$ is the integer defining $v=z_0^{k_0}$.  This limit depends only on our initially chosen $u$ and $z \in [u]$ (cf. Lemma \ref{lem:1}) and return words of $[v]$, but not   of the point $\omega$.
On the other hand, for $\nu$-almost all  $\omega'$,  Kingman's ergodic theorem implies
$$
    L(\omega')  = \lim_{n\to \infty} \frac{1}{n} \mathbb{E}_\nu \log \|A^{(n)}(x)\|.
$$
Thus, we have proved that the alternative (ii) of  Theorem \ref{Theorem-main-1} holds for any generic point $\omega$ with $A^{(n)}(\omega)\neq 0$ for all $n\in \N$.


\subsection{Proof of  lemma \ref{lemma-long-small}}\label{subsection-proof-lem-long-small}

Consider the average length of long return words defined by
\[
S^M_i=\frac{1}{\tau_i}\sum_{j=1}^{i} \left |\zeta_j\right|{\bf 1}_{R(v)\setminus F_v^M}(\zeta_j).
\]
For any fixed $i\ge 1$, the average $S_i^M$  is decreasing in $M$,  and so is  $\varlimsup_{i\to\infty} S_i^M$. We shall prove the conclusion of Lemma \ref{lemma-long-small} by contradiction.  Then suppose that there exists a $\delta>0$ such that 
$$
 \forall M\ge 1, \ \ \limsup_{i\to\infty} S_i^M\ge \delta.
 $$
 Let us fix an $M\ge 1$.  Let $\{i_\ell\}$ be a subsequence of
 integers  such that 
 \begin{equation}\label{eq:lemma-long-small 1} 
\lim_{\ell\to \infty} S_{i_\ell}^M\ge \delta.
 \end{equation}
For any $M\ge 1$ and $i\ge 1$ fixed,  consider the ``orbit measures along with long return words":
$$\nu_M^i:=\frac{1}{\tau_i}\sum_{j=1}^{i} {\bf 1}_{R(v)\setminus F_v^M}(\zeta_j)\left(\sum_{n=1}^{\tau_{j+1}-\tau_{j}-1}\delta_{\{\sigma^n(\zeta_j\zeta_{j+1}\cdots)\}}\right)$$
and the ``complementary orbit  measures" $\eta_M^i$ defined by 
$$\frac{1}{\tau_i}\sum_{n=0}^{\tau_i-1}\delta_{\sigma^n\omega} = \nu_M^i + \eta_M^i.$$
Recall that $\sigma^{\tau_j}\omega=\zeta_j\zeta_{j+1}\cdots$.
By the definition of $\tau_{j+1}$, for each $j\ge 1$ and $1\le n< \tau_{j+1}-\tau_j $, we have 
$\sigma^n(\zeta_j\zeta_{j+1}\cdots)\notin [v]$. Thus the measure $\nu_M^i$ doesn't charge the cylinder $[v]$, i.e., 
\begin{equation}\label{eq:lemma-long-small 2}
\forall M\ge 1, \forall i, \quad 
\nu_M^i([v])=0.
\end{equation}
Let $Y=\overline{\{\sigma^n\omega\}_{n\ge 1}}$ be the orbit closure of $\omega$.
In view of  \eqref{eq:lemma-long-small 1} , the fact $|\zeta_j|=\tau_{j+1}-\tau_j$ and the definition of $\nu_M^{i}$, we have
\begin{equation}\label{eq:lemma-long-small 3}
\forall M\ge 1, \quad \nu_M^{i_\ell}(Y)\ge \frac{\delta}{2}\ \  \textrm{ if } \ell \gg 1.
\end{equation}
Up to taking a subsequence of $\{i_\ell\}$, we can assume that 
$$\nu_M^{i_\ell}\rightharpoonup \nu_M^\infty \ \textrm{ and }\ \eta_M^{i_\ell} \rightharpoonup \eta_M^\infty$$
where $\nu_M^\infty$ and $\eta_M^\infty$ are some measures concentrated on $Y$. Note that $\nu_M^\infty$ and $\eta_M^\infty$ are not necessarily probability measures.  
Because of \eqref{eq:lemma-long-small 2}, we  have 
\begin{equation}\label{eq:lemma-long-small 4}
\forall M\ge 1, \quad \nu_M^\infty([v])=0.
\end{equation}
By \eqref{eq:lemma-long-small 3}, we have 
\begin{equation}\label{eq:lemma-long-small 4-1}
\forall M\ge 1, \quad 
\nu_M^{\infty}(Y)\ge \frac{\delta}{2}.
\end{equation}
Since $\omega$ is $\nu$-generic, the measure $\nu_M^i+\eta_M^i$ converges to $\nu$ in the weak star topology, as $i\to\infty$. Thus we have 
\begin{equation}\label{eq:lemma-long-small 4-2}
\forall M \ge 1, \quad \nu_M^\infty+\eta_M^\infty=\nu.
\end{equation}

In the following, we will study in more details the measure  $\nu_M^\infty$.  Associated to the $j$-th return word $\zeta_j$, we consider the orbit probability measure
$$
     \pi_j := \frac{1}{\tau_{j+1}-\tau_j}\sum_{n=1}^{\tau_{j+1}-\tau_j}\delta_{\{\sigma^n(\zeta_j\zeta_{j+1}\cdots)\}}.
$$
This is the orbit measure along with $\zeta_j$ (corresponding the finding of the next return word $\zeta_{j+1}$ or the returning to $[v]$).  
Let 
$$A_M=\left\{\pi_j: j\ge 1, \zeta_j\in R(v)\setminus F_v^M\right\},$$
the set of  orbit measures corresponding to  long return words. Notice that
$$
     \tau_{j+1}-\tau_j \ge M \quad {\rm when}\ \ \pi_j \in A_M.
$$  

\begin{lem} \label{lem:lim-mu}
Suppose that we have a sequence of measures $\mu_t \in A_{M_{t}}$, with $M_t \to \infty$ as $t\to \infty$, such that $\mu_t  \rightharpoonup \mu$. Then $\mu$ is $\sigma$-invariant. 
\end{lem}
Indeed, it is a direct consequence of the following known fact:  Given two sequences of integers
$\{p_k\}, \{q_k\}\subset\N$  with $q_k-p_k\to\infty $ as $k\to\infty$ and a sequence of points $\{x_k\} \in \{0, 1, \cdots, m-1\}^\mathbb{N}$. If 
$$\frac{1}{q_k-p_k}\sum_{j=p_k}^{q_k}\delta_{\sigma^j(x_k)}\rightharpoonup \lambda, \textrm{ as } k\to\infty,$$   
then the limit measure $\lambda$ is a $\sigma$-invariant measure.
\medskip

Continue our discussion. Let $\mathcal{A}$ be the weak-* closure of $\bigcup_M A_M$,   a compact set of measures. 
Let us rewrite $\nu_M^i$ as 
$$\nu_M^i=\frac{1}{\tau_i}\sum_{j=1}^{i}  (\tau_{j+1}-\tau_j) \pi_j {\bf 1}_{R(v)\setminus F_v^M}(\zeta_j).$$
Or equivalently 
$$\nu_M^i=\int_ {\mathcal{A}}\mu dQ_M^i(\mu)$$ 
where $Q_M^i$ is a discrete measure on $A_M (\subset \mathcal{A})$ with total mass not exceeding $1$.
Since $\nu_M^\infty$ is the weak-* limit of $\nu_M^{i_\ell}$,  
it can be written as 
\begin{equation}\label{eq:nu_M_infty}
\nu_M^\infty=\int_\mathcal{A} \mu dQ_M(\mu).
\end{equation}
for some  limit measure $Q_M$ of $Q_M^{i_\ell}$ (not necessarily probability measure) on the space $\overline{A}_M$. 
Let us explain this point.
Recall that equipped with the weak-* topology  the dual space $C^*(\Omega)$ is  locally compact and metrizable.  The set $\mathcal{A}$ is a compact subset in $C^*(\Omega)$. 
That $Q_M^{i_\ell}$ converges to $Q_M$ means 
$$
     \int_{\mathcal{A}} \varphi(\mu) d Q_M^{i_\ell}(\mu) \to \int_{\mathcal{A}} \varphi(\mu) d Q_M(\mu)
$$ 
for all weak-$*$ continuous function $\varphi$. In particular, as the function 
$\mu \to \int_\Omega fd\mu$ is weak-$*$ continuous, we have 
$$
    \int_{\mathcal{A}} \int_\Omega f d\mu d Q_M^i(\mu) \to \int_{\mathcal{A}} \int_\Omega f d\mu d Q_M(\mu). 
$$
So, the equality (\ref{eq:nu_M_infty}) reads as 
$$
\int_\Omega f d\nu_M^\infty =  \int_{\mathcal{A}} \int_\Omega f d\mu d Q_M(\mu)
$$ 
for all continuous function $f$ on $\Omega$.


Take a sequence of integers $\{M_t\}$ tending to $\infty$ such that
$$\nu_{M_t}^\infty\rightharpoonup\nu_\infty^\infty, \quad
\eta_{M_t}^\infty\rightharpoonup\eta_\infty^\infty , \quad Q_{M_t}\rightharpoonup Q_\infty \quad \textrm{ as } t\to \infty.$$
Then we  have 
\begin{equation}\label{eq:nu-inf}
\nu_\infty^\infty=\int_\mathcal{A} \mu d Q_\infty (\mu).
\end{equation}
 From \eqref{eq:lemma-long-small 4}, we get 
\begin{equation}\label{eq:lemma-long-small 5}
\nu_\infty^\infty([v])=0.
\end{equation}
From \eqref{eq:lemma-long-small 4-2}, we get
\begin{equation}\label{eq:lemma-long-small 6}
\nu_\infty^\infty+\eta_\infty^\infty=\nu.
\end{equation}

Since $Q_\infty $ is the weak limit of $Q_{M_t}$,  it holds that for each $\mu\in {\rm supp}(Q_\infty)$ and each $r>0$, we have $Q_{M_t}(B(\mu,r))>0$ for all large enough $t$, where 
$B(\mu, r)$ is the ball centered at $\mu$ of radius $r$.  It follows that for every $\mu\in {\rm supp}(Q_\infty)$, there is a sequence of probability measures $\{\mu_t\}$ with $\mu_t\in A_{M_t}$ such that   $\mu_t\rightharpoonup \mu$, as $t\to\infty$. 
By Lemma \ref{lem:lim-mu},
the measure $\mu$ is $\sigma$-invariant on $Y$.   Then $\nu_\infty^\infty$ is $\sigma$-invariant, by (\ref{eq:nu-inf}).

By \eqref{eq:lemma-long-small 4-2}, the invariant measure $\nu_\infty^\infty$ is absolutely continuous with respect to the measure $\nu$ which is assumed ergodic. It follows that
$\nu_\infty^\infty= c \nu$ for some constant $c$, which is not zero by (\ref{eq:lemma-long-small 4-1}).
Finally (\ref{eq:lemma-long-small 5}) contradicts $\nu([v])>0$.

\subsection{Proof of Theorem \ref{Theorem-main-1} when $\mu$ is atomic}\label{sect:atomic}

We assume that $A^{(n)}(\omega)\not=0$ for all $n\ge 0$.  Otherwise, $A^{(n)}(\omega)=0$ for large $n$ and the Lyapunov exponent is equal to $-\infty$. 

Let us first recall a simple fact about the spectral radius of matrix product. The spectral radius of a square matrix $C$ will be denoted $\rho(C)$.

\begin{lem}\label{lem:radius} Let $A$ and $B$ be two $d\times d$ matrices. We have $\rho(AB)=\rho(BA)$. 
\end{lem}
\begin{proof} It is trivial if $A$ or $B$ is the zero matrix. Otherwise, $\rho(AB)\le \rho(BA)$  follows from the Gelfand's formula 
$\rho(C) = \lim_n \|C^n\|^{1/n}$ and the estimate
$$
    \|(AB)^n\| \le \|A\| \|(BA)^{n-1}\|\|B\|;
$$
and then $\rho(AB)\ge \rho(BA)$ by symmetry. 
\end{proof}

Now let us start the proof. 
Recall that by Lemma~\ref{lemma-implication-condition-1},  there exists a word $u$ with  $|u|\ge \ell_0$ such that  
\[
 \nu([u])>0 \textrm{ and } b:=\min_{1\le i,j\le d} \min_{x\in [u]}A^{(\ell_0)}(x)_{i,j}>0.
\]
Assume that $\nu$ is the ergodic invariant measure supported by the orbit of a $p$-periodic point $x= (a_1a_2\cdots a_p)^\infty$ ($p\ge 1$ being an integer).  
The periodic cycle of $x$ is $\{x, \sigma x, \cdots, \sigma^{p-1} x\}=:Z $, namely
$$Z= \{(a_1a_2\cdots a_{p-1}a_p)^\infty, (a_2 a_3\cdots a_pa_1)^\infty,  \cdots,  (a_p a_1\cdots a_{p-2}a_{p-1})^\infty\}. $$
For any integer $\ell\ge 1$, consider the set of words
$$
    Z_\ell=\{(a_1a_2\cdots a_{p-1}a_p)^\ell, (a_2 a_3\cdots a_pa_1)^\ell, \cdots, (a_p a_1\cdots a_{p-2}a_{p-1})^\ell\}.
$$
The cylinders  $[z]$ with $z\in Z_\ell$ are neighborhoods of the $p$ points in the cycle.  Put them together to get a neighborhood
of the cycle:
$$
  \mathcal{Z}_\ell =\bigcup_{v\in Z_\ell}[v].
$$
Since $\nu$ is supported by the cycle and $\nu([u])>0$, $[u]$ contains at least one point in the cycle.  It follows  that if $\ell\ge 2$ and $p\ell>2|u|$, then for each $v\in Z_\ell$,  there exists $n\le \max(p,|u|)$ such that $u$ is a prefix of $\sigma^n(v)$.  Let $\ell_1$ be the least $\ell\ge 2$ such that $p\ell>2|u|$.

That $\omega$ is $\nu$-generic means that
\begin{equation}\label{eq:neglig}
  \forall \ell \ge 1, \quad    \lim_{n\to \infty}\frac{1}{n}\sum_{j=0}^{n-1}1_{\mathcal{Z}_\ell^c}(\sigma^j \omega)=0.
\end{equation}
That is to say, the patterns  different from those in $Z_\ell$
are negligible. 

By Lemma~\ref{lem:radius} we have 
$$
   \rho\left(A^{(p)}(z)\right)= \rho\left(A^{(p)}(z')\right) \textrm{  for all } z,z'\in Z.
$$
In the following we denote by {\color{blue}{$\rho>0$}} this common spectral radius of $A^{(p)}(z)$,  $z\in Z$. 
 Since the cycle is finite, by Gelfand's formula for spectral radius, for any
$\epsilon>0 $, there exist $\ell_2>\ell_1$ such that 
\begin{equation}\label{eq:rho}
 \forall \ell\ge \ell_2, \forall z \in Z,  \quad (\rho-\epsilon/2)^\ell \le\|A^{(\ell p)}(z)\| \le(\rho+\epsilon/2)^\ell.
\end{equation}
Since $A$ is continuous, we can perturb a little bit $z$, namely for large   $\ell_3\ge \ell_2$ and   for all $z\in Z$ and all $z'\in [z_0^{\ell_3-1}]$ (recall: $z_0^{\ell_3-1}$ denotes the prefix of $z$
of length $\ell_3$), we have
\begin{equation}\label{eq:rho-2}
 (\rho- \epsilon)^{\ell_2} \le\|A^{(\ell _2p)}(z')\| \le(\rho+\epsilon)^{\ell_2}
\end{equation}
and also
\begin{equation}\label{eq:rho-3}             
 \forall 1\le i,j\le d, \quad A(z')_{i,j}\ge (1-\epsilon)A(z)_{i,j}.
\end{equation}
In the following, we shall choose and  fix an $\ell_3=m\ell_2$ with $m\gg 1$.

Now,  we fix $\ell=k\ell_3$ with $k\gg 1$ and  we shall give a  decomposition of $\omega$, which is similar to that in (\ref{eq:decomp}).  Here we consider the return times to $\mathcal{Z}_\ell$, a small neighborhood of the periodic cycle. First, let $s_0$ be the least $s\ge 0$ such that $\sigma^s(\omega) \in \mathcal{Z}_\ell$. 
So, $\omega$ takes the form
$$
    \omega =\omega_0^{s_0-1}\omega_{s_0}^{s_0+\ell p-1 }\omega_{s_0+\ell p}^\infty =\mathbf{w}_0\mathbf{z}_0 \omega_{s_0+\ell p}^\infty,
$$ 
where $\mathbf{z}_0:=\omega_{s_0}^{s_0+\ell p-1 }$ is the first word from  $Z_\ell$ that we see in $\omega$ when we shift from the left, and  $\mathbf{w}_0:=\omega_0^{s_0-1}$ which is the empty word when $s_0=0$. 
Now we continue to look for the first word from  $Z_\ell$ that we see in $\omega_{s_0+\ell p}^\infty$. Let $s_1>s_0$ be the least $s\ge s_0+\ell p$ such that $\sigma^{s}(\omega) \in \mathcal{Z}_\ell$,  then let $s_2>s_1$  be the least $s\ge s_1+\ell p$ such that $\sigma^{s}(\omega) \in \mathcal{Z}_\ell$,  and so on.  By induction, we get a decomposition of $\omega$: 
$$
 \omega =\mathbf{w}_0\mathbf{z}_0  \mathbf{w}_1\mathbf{z}_1  \cdots \mathbf{w}_j\mathbf{z}_j \cdots 
$$
where for all $j\ge 1$ we have
\[
  \mathbf{w}_j= \omega_{s_{j-1} +\ell p}^{s_{j}-1}, \quad \mathbf{z}_j= \omega_{s_j}^{s_j+\ell p-1 }
\]
By construction, we have   $\mathbf{z}_j\in Z_\ell$ so that $\sigma^{s_j}(\omega)\in \mathcal{Z}_\ell$ for each $j\ge 0$.  But when shifted  over $\mathbf{w}_j$, the orbit of $\omega$
is outside $\mathcal{Z}_\ell$, namely $\sigma^i(\omega)\notin \mathcal{Z}_\ell$ for all $s_{j-1}+\ell p\le i\le s_{j}-1$.   This last fact, together with  \eqref{eq:neglig},  implies
\begin{equation}\label{eq:neglig-2}
 \ \lim_{r\to \infty}=\frac{|\mathbf{w}_0| + |\mathbf{w}_1|+\cdots +  |\mathbf{w}_r|}{s_r}=   \lim_{r\to \infty}\frac{s_0+\sum_{j=1}^{r}(s_{j}-s_{j-1}-\ell p)}{s_r} =0.
\end{equation}
This is equivalent to the fact $s_r \sim r p\ell$. In particular, we get 
\begin{equation}\label{eq:neglig-coro-1}
      \lim_{r\to \infty}\frac{s_{r+1}}{s_r}=1.
\end{equation}
According to this and Lemma \ref{lemma-basic-norm-estimate-1}, 
as in the case of continuous measure $\nu$,   
we only need to study 
$s_r^{-1}\log \|A^{(s_r)}(\omega)\|$ to show that it approaches $\frac{\log \rho}{p}$, when $\ell$ is sufficiently large.
Indeed, for any large  $n\ge 1$,  letting $r\ge 0$ be the integer such that $s_r\le n< s_{r+1}$,  by Lemma~\ref{lemma-basic-norm-estimate-1}, we have 
\[
\log \frac{\|A^{(n)}(\omega)\|}{\|A^{(s_r)}(\omega)\|}=O(n-s_r).
\]
Hence by \eqref{eq:neglig-coro-1},
\begin{equation}\label{eq:neglig-coro-2}
      \frac{1}{s_r}\log \|A^{(s_r)}(\omega)\| -\frac{1}{n}\log \|A^{(n)}(\omega)\| \to 0 \ \textrm{ as } n\to \infty.
\end{equation}

Recall that for each $j\ge 0$,  as we discussed at the beginning of the proof, there exists an integer $0\le n_j\le \max(p,|u|)$ such that $u$ is a prefix of $\sigma^{n_j}(\mathbf{z}_j)$,  so that
\[
\sigma^{n_j+s_j} (\omega)\in [u].
\]
That is to say, when the orbit $\sigma^{s_j}\omega$ falls into $\mathcal{Z}_\ell$ (i.e. when we see the word $\mathbf{z}_j$ from $Z_\ell$), it will  fall into $[u]$ if it goes forward  a small number $n_j$ steps.  
Thus,  by Lemma~\ref{lemma-implication-condition-1}, we have 
\begin{equation}\label{eq:pf-atomic-positivity-1}
\min_{1\le \alpha,\beta\le d}\left(A^{(|u|)}(\sigma^{n_j+s_j}\omega)\right)_{\alpha,\beta} \ge b >0.
\end{equation}
Since $\mathbf{z}_j=a^\ell$ for some word $a$ with $|a|=p$
and since $\ell=k\ell_3$,  we can write $\mathbf{z}_j = \mathbf{z}^k$ for  $\mathbf{z}=a^{\ell_3} \in Z_{\ell_3}$. 
Observe that $\sigma^{n_j}(\mathbf{z}_j)$ is of the form $w^{k-1}w'$ with $|w|=p \ell_3 $ and $|w'|=p\ell_3-n_j$ where $w$ has $u$ as prefix.
So,  $u$ is also a prefix of $\sigma^{(k-2)p\ell_3}\left(\sigma^{n_j}(\mathbf{z}_j)\right)$.  Thus
$\sigma^{n_j+s_j+(k-2)p\ell_3}\omega \in [u]$.
Again by Lemma~\ref{lemma-implication-condition-1}, we have 
\begin{equation}\label{eq:pf-atomic-positivity-2}
\min_{1\le \alpha,\beta\le d}\left(A^{(|u|)}(\sigma^{n_j+s_j+(k-2)p\ell_3}\omega)\right)_{\alpha,\beta} \ge b.
\end{equation}
In the following, for simplicity, we denote $L=(k-2)p\ell_3$, which is fixed.

With the help of 
\eqref{eq:pf-atomic-positivity-1} and \eqref{eq:pf-atomic-positivity-2} applied to $j=0$,  using  Lemma~\ref{lem:elem}  to 
break $\|A^{(s_r)}(\omega)\|$ twice, we get 
\[
\|A^{(s_r)}(\omega)\|= c^{\pm 2} \|A^{(n_0+s_0)}(\omega)\|\cdot
\|A^{(L)}(\sigma^{n_0+s_0}\omega)\|\cdot
\|A^{(n_0+s_0+L, s_r)}(\omega)\|,
\]
where $0<c<\infty$ is some constant depending only on $u$ and $A$.  In the same way, apply 
\eqref{eq:pf-atomic-positivity-1}, \eqref{eq:pf-atomic-positivity-2} with $j=1$ and  Lemma~\ref{lem:elem}  to get
\begin{eqnarray*}
\|A^{(n_0+s_0+L,s_r)}(\omega)\| & = & c^{\pm 2} \|A^{(n_0+s_0+L,n_1+s_1)}(\omega)\|\cdot
\|A^{(L)}(\sigma^{n_1+s_1}\omega)\| 
\cdot \|A^{(n_1+s_1+L,s_r)}(\omega)\|.
\end{eqnarray*}
Actually  for all $j=2, 3,\ldots,  r-1$, we also have 
\begin{eqnarray*}
\|A^{(n_{j-1}+s_{j-1}+L,s_r)}(\omega)\| & = & c^{\pm 2} \|A^{(n_{j-1}+s_{j-1}+L,n_{j}+s_{j})}(\omega)\|\cdot
\|A^{(L)}(\sigma^{n_{j}+s_{j}}\omega)\|  
 \cdot \|A^{(n_{j}+s_{j}+L,s_r)}(\omega)\|.
\end{eqnarray*}
Putting  the above estimates together and taking into account   \eqref{eq:neglig-2} and Lemma~\ref{lemma-basic-norm-estimate-1}, we get
\begin{equation}\label{eq:pf-atomic-quasi-multi-esti-1-1}
\|A^{(s_r)}(\omega)\|=c^{\pm 2r}e^{\pm   o(c_2s_r)}\prod_{j=0}^{r-1}\|A^{(L)}(\sigma^{n_j+s_j}\omega)\|.
\end{equation}

Now,  for each $0\le j\le r-1$,  we  write  
\begin{equation}\label{eq:pf-atomic-quasi-multi-esti-1-2-1}
A^{(L)}(\sigma^{n_j+s_j}\omega)=\prod_{i=0}^{k-3} A^{( p\ell_3)}(\sigma^{n_j+s_j+i\cdot p\ell_3} \omega).
\end{equation}
By the definition of $s_j$, the word $\omega_{s_j}^{s_j+kp\ell_3 -1}$ belongs to $Z_\ell=Z_{k\ell_3}$ and is $p$-periodic then $p\ell_3$-periodic.  Then there exists a periodic point $z\in Z$ such that for all $0\le i\le k-3$,
\begin{equation}\label{eq:pf-atomic-quasi-multi-esti-1-2-1-1}
  \forall q=0,1,\ldots,p\ell_3-1, \quad  \sigma^{n_j+s_j+i\cdot p\ell_3+q}\omega\in [\sigma^q(z)_0^{\ell_3-1}].
\end{equation}
Thus, by \eqref{eq:rho-3},  we have 
\begin{equation*}
\left(A^{( p\ell_3)}(\sigma^{n_j+s_j+i\cdot p\ell_3} \omega)\right)_{\alpha,\beta}\ge (1-\epsilon)^{p\ell_3}\left(A^{(p\ell_3)}(z)\right)_{\alpha,\beta},  \forall 1\le \alpha,\beta \le d
\end{equation*}
and consequently 
\begin{equation}\label{eq:pf-atomic-quasi-multi-esti-1-2-2}
\| A^{( p\ell_3)}(\sigma^{n_j+s_j+i\cdot p\ell_3} \omega)\|\ge (1-\epsilon)^{p\ell_3}\| A^{(p\ell_3)}(z)\|.
\end{equation}
Using  \eqref{eq:pf-atomic-quasi-multi-esti-1-2-1},  \eqref{eq:pf-atomic-quasi-multi-esti-1-2-2} and the first inequality in \eqref{eq:rho}, for all $ 0\le j\le r-1$ we then get
\begin{equation}\label{eq:pf-atomic-quasi-multi-esti-1-2}
\|A^{(L)}(\sigma^{n_j+s_j}\omega)\| \ge (1-\epsilon)^{L}\left(\rho-\epsilon/2\right)^{L/p}
\end{equation}
Recall that $L=(k-2)p\ell_3$.

On the other hand,   for every $0\le j\le r-1$, as $\ell_3=m\ell_2$ we have trivially
\begin{equation}\label{eq:pf-atomic-quasi-multi-esti-1-3}
 \|A^{(L)}(\sigma^{n_j+s_j}\omega)\| \le \prod_{i=0}^{m(k-2)-1} \| A^{( p\ell_2)}(\sigma^{n_j+s_j+i\cdot p\ell_2} \omega)\|.
\end{equation}
In view of  \eqref{eq:pf-atomic-quasi-multi-esti-1-2-1-1} ,  by the second inequality in  \eqref{eq:rho-2} we get 
\[ 
\| A^{( p\ell_2)}(\sigma^{n_j+s_j+i\cdot p\ell_2} \omega)\| \le (\rho+\epsilon)^{\ell_2}.
\]
Plug the above estimate into \eqref{eq:pf-atomic-quasi-multi-esti-1-3},   for all $ 0\le j\le r-1$ we obtain
\begin{equation}\label{eq:pf-atomic-quasi-multi-esti-1-4}
\|A^{(L)}(\sigma^{n_j+s_j}\omega)\| \le (\rho+\epsilon)^{\ell_2  m(k-2)} =(\rho+\epsilon)^L.
\end{equation}

Combining \eqref{eq:pf-atomic-quasi-multi-esti-1-1},  \eqref{eq:pf-atomic-quasi-multi-esti-1-2} and \eqref{eq:pf-atomic-quasi-multi-esti-1-4},  we get
\[
o(s_r)+r L[\log (1-\epsilon)+\log(\rho-\epsilon/2)]
\le \log \|A^{(s_{r})}(\omega)\| 
\le rL \log (\rho+2\epsilon)+o(s_r). 
\]
Thus, as $s_r \sim r p\ell = rk p\ell_3\sim r L$ (cf. \ref{eq:neglig-2}), we have
\[
\liminf_{r\to\infty} \frac{\log \|A^{(s_r)}(\omega)\|}{s_r} 
\ge \frac{\log \rho}{p}-o_\epsilon(1)-o_k(1),
\]
where $o_\epsilon(1)\to 0$ as $\epsilon\to 0$ and $o_k(1)\to 0$ as $k\to \infty$.  Similarly, from  (\ref{eq:pf-atomic-quasi-multi-esti-1-4}) we get  
\[
\limsup_{r\to\infty} \frac{\log \|A^{(s_r)}(\omega)\|}{s_r} 
\le \frac{\log \rho}{p}+o_\epsilon(1)+o_k(1).
\]
First take  $\epsilon>0$ very small (by choosing $\ell_3\gg1$) and then take $k$ very large, so that the above liminf and limsup can be made very close
to $\frac{\log \rho}{p}$.  Finally, 
from \eqref{eq:neglig-coro-1},  we  get 
\[
\lim_{n\to\infty}\frac{\log \|A^{(n)}(\omega)\|}{n}=\frac{\log \rho}{p}.
\]

\section{Remarks and Counter-examples}\label{sect:remarks}
We would like to discuss the positivity conditions (\ref{eq:claim1-cond-1}) and  (\ref{eq:claim1-cond-2}) in our main Theorem  \ref{Theorem-main-1}, which ensures the existence of Lyapunov exponent (\ref{eq:Liap0}).
 We start with recalling the Birkhoff contraction coefficient  in terms of Hilbert projective metric on allowable non-negative matrices (cf.  \cite{Seneta}), and comparing the condition (\ref{eq:claim1-cond-1})
 and the contractivity condition which is usually used for studying non-negative matrices (cf. \cite{Hajnal1976}). 
\medskip 

\subsection{Contraction in terms of Hilbert projective metric}
A nonnegative matrix $B=(b_{i, j})\ge 0$ (i.e. the entries $b_{i, j}\ge 0$) is characterized by the fact that it preserves the cone
 $K$ in the sense that $BK\subset K$ where 
 $$
     K=\left\{(x_1, \cdots, x_d)\in \mathbb{R}^d: x_i\ge 0\right\}.
 $$ 
 We write $ \mathring{K} $ for the interior of $K$ and $K^+$ for $K\setminus\{0\}$.  A row-allowable non-negative matrix $B$ defines a map $B:\mathring{K} \to \mathring{K}$, which  is contractive respect to the 
 {\em  Hilbert projective metric} defined by 
 $$
      d_H(x, y) = \log \frac{\max_i (x_i/y_i)}{\min_i (x_i/y_i)} \qquad {\rm for}\  x=(x_i) \in \mathring{K},\  y=(y_i) \in \mathring{K}.
 $$
 The function $d_H$ has all the properties of a metric, with one exception that  $d_H(x, y) =0$ if and only if $x=\lambda y$ for some $\lambda >0$. 
 The contraction property  of a row-allowable $B\ge 0$ is stated  as $$
      d_H(Bx, By)\le d_H(x, y)\quad {\rm  for\  }   x, y \in \mathring{K}.
      $$   
      (See Lemma 3.1 in \cite{Seneta}).
      This contraction is a very important fact. The {\em Birkhoff contraction coefficient} of $B$ is then defined by
 $$
      \tau(B) =\sup_{x,y\in \mathring{K}; x\not= \lambda y} \frac{d_H(Bx, By)}{d_H(x, y)}.
 $$
 It is easy to see that 
     $$0\le \tau(B)\le 1, \qquad \tau(B) =\tau(B'), \qquad \tau(B_1B_2)\le \tau(B_1)\tau(B_2)$$ 
     for any row-allowable matrices $B\ge 0, B_1\ge 0, B_2\ge 0$ ($B'$ denoting the transposed matrix of $B$). 
 Another important fact is that the Birkhoff contraction coefficient of an allowable matrix $B=(b_{i, j})$ can be expressed by cross-ratios of the entries of $B$,
 namely 
 \begin{equation}\label{eq:formula}
     \tau(B) =\frac{1-\sqrt{\phi(B)}}{1+\sqrt{\phi(B)}}
 \end{equation}
 where $\phi(B) =0$ when $B$ has at least one zero entry, and when $B>0$ we have
 $$
     \phi(B) = \min_{i, j, k, s} \frac{b_{i,k} b_{j,s}}{b_{j,k}b_{i,s}}.
 $$
 See the section 3.4 of \cite{Seneta} for a proof of (\ref{eq:formula}). 
 We say that a sequence of allowable non-negative matrices $\{B_{n}\}$ is  {\em forward contractive} if
      \begin{equation}\label{eq:F-ergodicity}
     \lim_{n \to \infty} \tau(B_1B_{2}\cdots B_{n}) =0.
     \end{equation}
     The terminology "weak ergodicity"
 is used in \cite{Seneta}, pp. 84. We prefer to reserve the "ergodicity" for its proper sense in  the ergodic theory. 
     The backward contractive property is similarly defined by  $\tau(B_{n}B_{n-1}\cdots B_{1}) \to 0$. 
      A sequence of allowable non-negative matrices is said to be {\em contractive} if it is both forward contractive and backward contractive.
      \medskip

\subsection{Positivity condition \eqref{eq:claim1-cond-1}  in Theorem  \ref{Theorem-main-1} compared with contractivity.}  Assume that $A_0, A_1, \cdots, A_{m-1}$ are row-allowable non-negative matrices. 
Here we consider the function $A(\cdot)$ depending on the first coordinate defined by $A(x) = A_{x_0}$ for $x=(x_{n})_{n\ge 0}$.  The condition \eqref{eq:claim1-cond-1} implies that for any 
$\nu$-generic point $\omega$,
the sequence $\{A_{\omega_n}\}$ is forward contractive: 
$$
\tau(A_{\omega_0}A_{\omega_1}\cdots A_{\omega_{n-1}})\to 0, \quad {\rm  as}  \ n\to \infty.
$$ Indeed, by the assumption there exists a word  $u$ such that  $\nu([u])>0$ and $\tau(A_u)<1$.  Since $\omega$ is $\nu$-generic, the word $u$ appears infinitely many times in the sequence $\omega$.  Using the fact $\tau(AB)\le \tau(A)\tau(B)$,  we deduce the forward contractivity of $\{A_{\omega_n}\}$.  But, in general, the converse is not true. See the following example.

\begin{exm}
Assume $m=2$, $A_0 =\begin{pmatrix} 1&0\\0&1\end{pmatrix}$ and $A_1=\begin{pmatrix} 1&1\\1&1\end{pmatrix}$,  $\nu=\delta_{0^\infty}$ which is ergodic  and $\omega=10^\infty$ which is $\nu$-generic
(i.e. every pattern $0^\ell$ has frequency $1$). We have $\tau(A_{\omega_0}A_{\omega_1}\cdots A_{\omega_{n-1}})= \tau(A_1)=0$ for all $n\ge 1$, but 
the condition (\ref{eq:claim1-cond-1}) is not satisfied. Actually there are uncountably many such $\nu$-generic points. For example, $\omega=(\omega_n)_{n\ge 0}$ defined by $\omega_n =0$ for $n\not=2^k$ ($k\ge 0$) and 
$\omega_{2^k}=1$ for at least one $k \ge 0$.
\end{exm}

However,  if  $\omega$ belongs to the support of $\nu$, 
 then the condition \eqref{eq:claim1-cond-1} is equivalent to that  $\{A_{\omega_n}\}$ is forward contractive, because any word $\omega_0\omega_1\cdots \omega_{n-1}$ appears infinitely often in $\omega$.
 
 Hajnal \cite{Hajnal1976} proved that if non-negative matrices $B_1, B_2, \cdots, B_{n}, \cdots$ are allowable, then the   forward contractivity of the sequence 
 $\{B_n\}$  means the product $B_1B_2\cdots B_n$
  tends to row-proportionality as $n\to \infty$  (cf. Theorem 1 in \cite{Hajnal1976}).  So, the contractivity describes well an asymptotical property of the products 
 $B_1B_2\cdots B_n$.  But the existence of Lyapunov exponent $\lim_{n\to \infty} n^{-1}\|B_1B_2\cdots B_n\|$, which is another asymptotical property, can not be ensured by the contractivity, as Example \ref{exm:nolimit} below shows. 

 \subsection{Necessity of the positivity \eqref{eq:claim1-cond-1}  in Theorem  \ref{Theorem-main-1} }
The positivity condition \eqref{eq:claim1-cond-1}  in Theorem  \ref{Theorem-main-1} can be dropped. Walters' counter-example shows that. 
Here we present a simpler counter-example on  $\{0,1,2,3\}^\N$ with $A(\cdot)$ depending only on the first coordinate,  which will satisfy:
\begin{itemize}
\item[(1)]  $\nu$ is an ergodic measure on $(\{0,1,2,3\}^\N,\sigma)$ (which can be continuous or atomic),
\item[(2)] $\omega\in  \{0,1,2,3\}^\N$ is a generic point of $\nu$,
\item[(3)] $A_0, A_1, A_2,A_3$ are $2\times 2$ non-negative and allowable matrices,
\item[(4)] 
$\lim_{n\to \infty}\tau(A_{\omega_1}\cdots A_{\omega_n}) =\lim_{n\to \infty} \tau(A_{\omega_n}A_{\omega_{n-1}}\cdots A_{\omega_1})=0$,
\item[(5)] The limit $\lim_{n\to\infty}\frac{1}{n}\log \|A_{\omega_1}\cdots A_{\omega_n}\|$ doesn't exist.
\end{itemize}
Here is such a example.

\begin{exm}\label{exm:nolimit}
Let $\nu$ be an ergodic measure (continuous or atomic) whose support is included in the subset $\{0,1\}^\N\subset \Omega$. Let $x\in \{0,1\}^\N$  be a $\nu$-generic point. 
Define 
$$A_0=A_1=\left( \begin{array}{ccc} 10 & 0\\ 0 & \frac{1}{10} \end{array} \right),
\quad
A_2=\left( \begin{array}{ccc} 0 & 1\\ 1 & 0 \end{array} \right), \quad
A_3=\left( \begin{array}{ccc} 1 & 1\\ 1 & 1 \end{array} \right).
$$
Define
$$ \omega=\omega_0 u_1v_1u_2v_2\cdots u_nv_n \cdots
$$
 where $\omega_0=3$ and 
$$
u_n =v_n= \left\{ \begin{array}{ll}
x_1\cdots x_n & \textrm{if $2^{2^{2k}}\le n<2^{2^{2k+1}},  k\in \N$}\\
x_1\cdots x_n2 & \textrm{if $2^{2^{2k+1}}\le n<2^{2^{2k+2}},  k\in \N$.}
\end{array} \right.
$$
Since $\tau(A_3)=0$,  we have $\tau(A_{\omega_0}A_{\omega_1}\cdots A_{\omega_{n-1}})= 0$ and $\tau(A_{\omega_{n-1}}\cdots A_{\omega_1}A_{\omega_0})= 0$.
Observe that 
$$
(A_{x_1}\cdots A_{x_n})^2=\left( \begin{array}{ccc} 10^{2n} & 0\\ 0 &  \frac{1}{10^{2n}} \end{array} \right), \quad
(A_{x_1}\cdots A_{x_n}A_2)^2=\left( \begin{array}{ccc} 1& 0\\ 0 &  1 \end{array} \right).$$  Using these facts, it is easy to check that $\frac{1}{n}\log \|A_{\omega_0}A_{\omega_1}\cdots A_{\omega_{n-1}}\|$ does not converge as $n\to\infty$.
\end{exm}

This example also shows that the contractivity can not replace the positivity \eqref{eq:claim1-cond-1}.

\subsection{Necessity of the positivity condition \eqref{eq:claim1-cond-2}  in Theorem  \ref{Theorem-main-1}}
We  construct now an example that satisfies all the conditions of Theorem  \ref{Theorem-main-1} except \eqref{eq:claim1-cond-2},  while there exists a $\nu$-generic point $\omega\in \Omega:=\{0,1\}^\N$ such that the Lyapunov exponent $L(\omega)$ does not exist.  
Recall that the usual  metric on $\Omega$ is defined by $d(x,y)=2^{-\min(k:x_k\neq y_k)}$.

\begin{exm}\label{exm:f(x)} Let  $\Omega=\{0,1\}^\N$. Let $f: \Omega\to \R$ be the continuous function defined as follows
\[f(0^\infty)=0 \ \textrm{ and } f(x)=\exp\left(-d(0^\infty,x)^{-1}\right) \textrm{ if } x\neq 0^\infty. \]   
Then we define the continuous map $A: \Omega \to M_2(\R)$ by
$$
A(x)\mapsto 
\left( \begin{array}{ccc}
f(x) & f(x)   \\
f(x)   & f(x) 
\end{array} \right).
$$
Consider the ergodic measure  $\nu=\frac{1}{2}\delta_{(01)^\infty}+\frac{1}{2}\delta_{(10)^\infty}$ supported by the 2-periodic cycle $\{(01)^\infty,  (10)^\infty\}$ and the $\nu$-generic point  
$$
\omega=(01)^{n_1}0^{n_1'}\cdots (01)^{n_i}0^{n_i'}\cdots $$
 where $n_i=2^{2^{2^i}}$ and $n_i'=\log_2 n_i$.
Then,  it is not hard to check that 
\begin{itemize}
\item[(1)] the point $\omega$ is  $\nu$-generic;
\item[(2)] all the assumptions  except \eqref{eq:claim1-cond-2} in Theorem  \ref{Theorem-main-1} are satisfied;
\item[(3)] we have 
\[\liminf_{n\to\infty}\frac{1}{n} \log \|A^{(n)}(\omega)\|=-\infty \quad \textrm{ but} \quad \limsup_{n\to\infty}\frac{1}{n} \log \|A^{(n)}(\omega)\|> -\infty
\]
\end{itemize}
\end{exm}

For the above generic point $\omega$, the quasi additivity in Theorem \ref{thm:main1-1} is not satisfied for $\varphi_n(x) =\log \|A^{(n)}(x)\|$. But for the periodic point $(01)^\infty$, 
the quasi additivity  is satisfied. 

\subsection{Necessity of ergodicity in Theorem  \ref{Theorem-main-1}}

First we see that 
Lemma \ref{lemma-long-small}, which is a key  for the proof of  Theorem  \ref{Theorem-main-1}, does not hold when $\nu$ is only assumed invariant.

\begin{exm} \label{exm:gap}
Let  $\Omega=\{1,2,3,4\}^\N$. Let $\nu_1$ and  $\nu_2$ be two continuous  ergodic measures respectively supported by $\{1,2\}^\N$ and $\{3,4\}^\N$.  Consider the non-ergodic invariant measure 
 $\nu=\frac{1}{2}(\nu_1+\nu_2)$.
Take a $\nu_1$-generic point  $x\in \{1,2\}^\N$ and a $\nu_2$-generic point $y\in \{3,4\}^\N$.  Then construct the point
$$
z=u_1v_1u_2v_2\cdots u_nv_n\cdots, 
$$
where $u_n=x_0x_1\cdots x_{n-1}$ and $v_n=y_0y_1\cdots y_{n-1}$. We can check that
\begin{itemize}
\item[(1)] $z$ is a generic point for $\nu$.   
\item[(2)] The conclusion of Lemma \ref{lemma-long-small} does not hold true for $\omega=z$.
\end{itemize}
\end{exm}


In general, Theorem \ref{Theorem-main-1} does not hold if $\nu$ is only assumed invariant. Here is a counter-example.
\begin{exm}
Let $\Omega=\{1,2,3,4\}^\N$.  Take take $a=10,b=1/10$ and define
$$A_1=\left( \begin{array}{ccc} a & 0\\ 0 & b \end{array} \right), \quad
A_2=\left( \begin{array}{ccc} \frac{1}{b} & 0\\ 0 & \frac{1}{a} \end{array} \right),\quad
A_3=\left( \begin{array}{ccc} 1 & 1\\ 1 & 1 \end{array} \right),\quad
A_4=\left( \begin{array}{ccc} 0 & 1\\ 1 & 0 \end{array} \right). 
$$
It is clear that 
$$A_1^nA_2^n=\left( \begin{array}{ccc} \frac{a^n}{b^n} & 0\\ 0 &  \frac{b^n}{a^n} \end{array} \right)=
\left( \begin{array}{ccc} 100^n & 0\\ 0 &  \frac{1}{100^n} \end{array} \right), \quad 
A_1^nA_4A_2^nA_4=\left( \begin{array}{ccc} 1& 0\\ 0 &  1 \end{array} \right).
$$ 
Consider the non-ergodic invariant measure 
$
\nu=\frac{1}{3} (\delta_{1^\infty} + \delta_{2^\infty} + \delta_{3^\infty})
$
and the point 
$$\omega=u_1v_1w_1 u_2v_2w_2\cdots u_nv_nw_n \cdots$$
 where the words $u_n, v_n$ and $w_n$ are defined as follows
 \begin{eqnarray*}
w_n&=&3^n,\\
u_n &=& \left\{ \begin{array}{ll}
1^n & \textrm{if} \ 2^{2^{2k}}\le n<2^{2^{2k+1}},  k\in \N\\
1^n4 & \textrm{if}\  2^{2^{2k+1}}\le n<2^{2^{2k+2}},  k\in \N,
\end{array} \right.\\
v_n &=& \left\{ \begin{array}{ll}
2^n & \textrm{if}\  2^{2^{2k}}\le n<2^{2^{2k+1}},,  k\in \N\\
2^n4 & \textrm{if}\ 2^{2^{2k+1}}\le n<2^{2^{2k+2}},  k\in \N.
\end{array} \right.
\end{eqnarray*}
It can be checked that 
\begin{itemize}
\item[(1)] all the assumptions of Theorem \ref{Theorem-main-1} are satisfied, except that $\nu$ is ergodic;
\item[(2)]
$\frac{1}{n}\log \|A_{\omega_1}\cdots A_{\omega_n}\|$ does not converge as $n\to\infty$.
\end{itemize}
\end{exm}

\subsection{Two typical classes of generic points}
Here are two  typical examples to which Theorem \ref{Theorem-main-1} will apply.  We  are given two non-negative matrices $A_0$ and $A_1$  ($m=2$). The first example is    the Thue-Morse sequence defined by 
	$$\omega_n=\frac{1}{2}(1+ (-1)^{s_2(n)}),$$ where 
	$s_2(n)= \sum \epsilon_k$ is the sum of digits of $n$ in its dyadic expansion $n=\sum_k \epsilon_k 2^k$ ($\epsilon_k =0$ or $1$). More generally our result applies to all uniquely ergodic subshifts, thus including all primitive substitutive subshifts, Sturmian subshifts and linearly recurrent subshifts. 
	 The second example is the square of the M\"{o}bius function $\omega_n = \mu(n)^2$. 
	Recall that the  M\"{o}bius function $\mu$ can be defined by  the inverse of the Riemann zeta function:
	$$
	    \frac{1}{\zeta(s)} = \sum_{n=1}^\infty \frac{\mu(n)}{n^s}.
	$$
	Notice that the Thue-Morse sequence is a typical example of primitive substitutive sequences and it generates a subshift which is minimal and uniquely ergodic (cf. \cite{Queffelec}), but the square of the M\"{o}bius function 
	is of other kind and generates a subshift which is neither
	minimal nor uniquely ergodic, however the patterns contained in the  square of the M\"{o}bius function is well described by the so-called Mirsky measure which is ergodic and of zero entropy
	and $\mu(n)^2$ is generic for the Mirsky measure
	(cf. \cite{ALR2015, Peckner2015, Sarnak2012}). This genericity relative to the Mirsky measure ensures the applicability of Theorem \ref{Theorem-main-1}  to $\mu(n)^2$. 
	The square of the M\"{o}bius function is an special example of the characteristic functions of  $\mathcal{B}$-free integers  (cf. \cite{ALR2015,KLW2015}), to which Theorem \ref{Theorem-main-1} applies too. 

\subsection{Other remarks} Let us make the following remarks.
\medskip

R1. In Theorem \ref{thm:main1}, 
the measure $\nu$ is only assumed invariant, not necessarily ergodic. But we assume the strong positivity
(i.e. $A(y)>0$ for all $y$). 
This condition is technical and should be weakened. We are succeeded in controlling the gaps between the occurrences of a fixed  positive partial product in the shift dynamics for ergodic measure (cf. Lemma \ref{lemma-long-small}). Observe that the ergodicity can not be dropped
for controlling the gaps (cf. Example \ref{exm:gap} above).   Lemma \ref{lemma-long-small} can be generalized to other dynamics, but efforts are needed to generalize Theorem  \ref{Theorem-main-1}  to other dynamics.  
\medskip


R2. 
Let $\mathcal{L}(\omega)$ be the set of all finite words contained in $\omega$ (called the language of $\omega$). Theorem \ref{thm:main1-1} has the following immediate corollary. 
\begin{cor}\label{cor:Gen}
Suppose we have a function $\varphi: \mathcal{L}(\omega) \to \mathbb{R}$ satisfying the conditions 
\begin{itemize}
\item [(i)] $\varphi(uv) = \varphi(u) + \varphi(v) +o(|u|\wedge |v|)$ for all $u,v \in  \mathcal{L}(\omega) $ such that $uv\in  \mathcal{L}(\omega)$;
\item [(ii)] $\omega$ is a $\nu$-generic point for some shift invariant measure $\nu$.
\end{itemize}
Then the limit $\lim_{n\to \infty} \frac{\varphi(\omega_1\omega_2\cdots \omega_n)}{n}$ exists.
\end{cor}

Given a subshift $(\Omega, \sigma)$ over a finite alphabet. Let $\mathcal{W}(\Omega)$ be the set of finite words associated to $\Omega$. 
A function $\varphi: \mathcal{W}(\Omega) \to \mathbb{R}$ is said to be almost additive if there exist a constant $D>0$  and a non-increasing function $c: [0, \infty)\to [0, \infty) $
with $c(r) =o(1)$ as $r\to \infty$ such that\\
\indent \  (A1) \ $|\varphi(v) - \sum_{j=1}^n \varphi(v_j)|\le \sum_{j=1}^n c(|v_j|) |v_j|$ for $v=v_1v_2\cdots v_n\in \mathcal{W}(\Omega)$;\\
\indent \  (A2) $ |\varphi(v)|\le D|v|$ for $v\in \mathcal{W}(\Omega)$.\\
 Lenz proved  that the unique ergodicity of $(\Omega, \sigma)$ implies 
the existence of $\lim_{|w|\to \infty} \frac{\varphi(w)}{|w|}$ for every almost additive function $\varphi$ and that the inverse is also true (cf. Theorem 1 in \cite{Lenz2002}).

The quasi additivity  (i) in Corollary  \ref{cor:Gen} is stronger than the almost additivity. But Corollary \ref{cor:Gen} does not require the unique ergodicity of the subshift
and the condition (i) is local, namely only involves the orbit of the single point $\omega$. 

Recall that a function $\varphi: \mathcal{W}(\Omega) \to \mathbb{R}$ is said to be subadditive if $\varphi(uv)\le \varphi(u) +\varphi(v)$. 
Under the assumption that $(\Omega, \sigma)$ is a minimal subshift, Lenz also proved that $\lim_{|w|\to \infty} \frac{\varphi(w)}{|w|}$ exists for every subadditive function  $\varphi$
if and only if the subshift $(\Omega, \sigma)$ has the so-called uniform positivity of quasiweight (PQ), namely there exists a constant $C>0$ such that
$$
    \forall v \in \mathcal{W}(\Omega), \quad   \varliminf_{|w|\to \infty}\frac{N^*_v(w)}{|w|}\ge \frac{C}{|v|}
$$ 
where $N^*_v(w)$ is the maximal number of disjoint copies of $v$ in $w$ (cf. Theorem 2 in \cite{Lenz2002}). Notice that $N^*_v(w)$ is bounded from above by $N_v(w)$, the number of appearances of $v$ in $w$.
Subshifts  $(\Omega, \sigma)$ having the (PQ) property are uniquely ergodic and linear recurrent subshifts have the the (PQ) property.   
\medskip

R3. The formula (\ref{eq:Liap-Formula}) provides a kind of approximation for the  Lyapunov exponent. It is valid for independent and identically distributed non-negative random matrices,
under some positivity condition.   
In this case, the formula (\ref{eq:Liap-Formula}) is not practical because any word is a return word; but the transfer operator is successfully used in the computation of Lyapunov exponent by Pollicott \cite{Pollicott2010}. 
However, for simple systems like minimal systems, there is a finite number of return word and the formula (\ref{eq:Liap-Formula}) would have its merits. 
\medskip

R4. Theorem \ref{thm:main1}, as well as Theorem   \ref{Theorem-main-1}, can be used to construct Ising models for which the interactions are not constant, but form a generic point of shift-invariant measure. 
The free energy and the Gibbs measure are then well defined. 
\medskip

R5. Theorem \ref{thm:main1}  can be generalized to vector bundles. 
Let $X$ be a compact metric space and $T:X \to X$ be a homeomorphism and  $\pi: E \to X$ be a $d$-dimensional vector bundle over $X$ ($d\ge 1$). Let
$E_x = \pi^{-1}(x)$, which is a  $d$-dimensional vector space. Fix a Riemannian metric on $E$ (all metric are equivalent). 
Suppose that $A: E \to E$ is a vector bundle automorphism of $E$ covering $T$. Therefore $A$ maps linearly $E_x$ to $E_{Tx}$. We denote this linear map by $A_x$. 
The cocycle $(A^n)_x: E_x \to E_{T^nx}$ is defined by the linear map
$$
 (A^n)_x = A_{T^{n-1}x} \circ \cdots \circ  A_{Tx} \circ A_x.
$$ 
Its norm $\|(A^n)_x\|$ is defined as the norm of linear operator using the norms on $E_x$ and $E_{T^nx}$ given by the Riemannian metric.

Let $\mathcal{P}(E,T)$ be the set of all vector bundle automorphisms of $E$ covering $T$ and satisfying the following condition: for any $x\in X$
there exists a proper cone $C_x$ in $E_x$ such that $C_x$ varies continuously with $x$ and 
$$
     A_x(C_x\setminus\{0\}) \subset {\rm int} (C_{Tx}).
$$

Here by a cone we mean a subset $C$ of a vector space such that $C+C\subset C$ and $\alpha C \subset C$ for all $\alpha >0$. A cone is closed if it is a closed
subset of the (normed) vector space and a closed cone $C$ is proper if $C \cap (-C)=\{0\}$. 
By ${\rm int} (C)$ we mean the interior of $C$. The above continuity condition means that in the sphere bundle of $E$ (which is compact), the set determined by 
$C_x$ varies continuously in the space of all compact sets of the sphere bundle equipped with the Hausdorff metric.

\begin{thm}
 Let $T:X\to X$ be a homeomorphism of a compact metric space and  $\pi: E\to X$ be a $d$-dimensional vector bundle over $X$. 
 Suppose $A\in \mathcal{P}(E, T)$ and $y$ is a generic point for some $T$-invariant measure $\mu$, then the following limit exists
 \begin{equation}\label{eq:w2}
       \lim_{n\to \infty} \frac{1}{n} \log \|(A^n)_y\|.
 \end{equation}
\end{thm}

This result is actually proved by Walters \cite{Walters1986}. We just need to point out that at the end of his proof, instead of using Oxtoby's theorem, one uses the 
$\mu$-genericity of the point $y$. Walters' proof is based on the following result due to Ruelle (\cite{Ruelle1979}): there exists a nowhere-zero section $s(x)$ of $E$
and a strictly positive continuous function $a(x)$ on $X$ such that 
$$
      \forall x\in X, \quad A s(x) = a(x) s(Tx),
$$ 
and that  there exists a sub-bundle $\{W_x\}$ of $E$ with the property
$$
    \forall x\in X, \quad   AW_x =W_{Tx}, \quad E_x = F_x\oplus W_x
$$ 
where $F_x$ is  the one-dimensional space generated by $s(x)$, and there exist constants $0<\alpha<1$ and $K>0$ such that
$$
     \|A^n w\|\le K\alpha^n \|A^n u\|
$$
for all $w\in W_x$ with $\|w\|=1$, for all $u\in F_x$ with $\|u\|=1$, for all $n\ge 0$ and for all $x\in X$. Walters' proof also shows that the limit in \eqref{eq:w2} is equal to
$$
     \int_X \log a(x) d\mu(x).
$$

  \end{document}